\documentclass[a4,11pt]{article}
\usepackage[latin1]{inputenc}
\usepackage{amsmath,amsbsy,amsfonts,amssymb}
\usepackage{epsfig}
\usepackage{amsmath}
\usepackage{amssymb}
\usepackage{graphicx}

\newtheorem{theo}{Theorem}[section]

\newtheorem{lema}[theo]{Lemma}
\newtheorem{coro}{Corollary}[theo]

\newcommand{\fim}{{\hfill $\rule{2.0mm}{2.0mm}$}}
\newcommand{\proof}{\noindent{\bf Proof.\/ }}

\def \dsp {\displaystyle}

\def\XXint#1#2#3{{\setbox0=\hbox{$#1{#2#3}{\int}$}
     \vcenter{\hbox{$#2#3$}}\kern-.5\wd0}}

\def\downbar#1{
\setbox10=\hbox{$#1$}
            \dimen10=\ht10 \advance\dimen10 by 2.5pt
            \ifdim \dimen10<15pt 
               \advance\dimen10 by -0.5pt
               \dimen11=\dimen10
               \advance\dimen10 by 2.5pt
               \lower \dimen11
            \else \lower \ht10 \fi
            \hbox {\hskip 1.5pt \vrule height \dimen10 depth \dp10}}
\def\upbar#1{
\setbox10=\hbox{$#1$}
            \dimen10=\ht10 \advance\dimen10 by \dp10 \advance\dimen10 by 2.5pt
            \ifdim \dimen10<15pt 
               \advance\dimen10 by 2pt \fi
            \raise 2.5pt \hbox {\hskip -1.5pt \vrule height \dimen10}}

\begin{document}

\title{\Large{ Orthogonal polynomials on the unit circle:  Verblunsky coefficients with some restrictions  imposed on a pair of related real sequences}\thanks{The first and third authors are supported by funds from FAPESP (2014/22571-2) and CNPq (475502/2013-2, 305073/2014-1, 305208/2015-2) of Brazil. The second and fourth authors are supported by grants  from CAPES of Brazil.}}

\author
{
 {Cleonice F. Bracciali$^a$, Jairo S. Silva$^{b,}$\thanks{jairo.santos@ufma.br (corresponding author).}, A. Sri Ranga$^a$, Daniel O. Veronese$^c$}
  \\[0.5ex]
 {\small $^a$Departamento de Matem\'{a}tica Aplicada, IBILCE, } \\
 {\small UNESP - Universidade Estadual Paulista, } 
 {\small 15054-000, S\~{a}o Jos\'{e} do Rio Preto, SP, Brazil. } \\[0.5ex]
{\small $^b$Depto de Matem\'{a}tica, Universidade Federal do Maranh\~{a}o, 65080-805, S\~{a}o Lu\'{\i}s, MA, Brazil}\\
{\small and Pós-Graduação em Matem\'{a}tica, IBILCE, } \\
 {\small UNESP - Universidade Estadual Paulista, } 
 {\small 15054-000, S\~{a}o Jos\'{e} do Rio Preto, SP, Brazil. }\\[0.5ex]
 {\small $^c$ICTE, Universidade Federal do Tri\^{a}ngulo Mineiro,
38064-200, Uberaba, MG. }\\[0ex]
}

\date{ }

\maketitle

\thispagestyle{empty}

\begin{abstract}
It was shown recently that associated with a pair of real sequences $\{\{c_{n}\}_{n=1}^{\infty}, \linebreak \{d_{n}\}_{n=1}^{\infty}\}$, with $\{d_{n}\}_{n=1}^{\infty}$  a positive chain sequence, there exists a unique nontrivial probability measure $\mu$ on the unit circle. The Verblunsky coefficients $\{\alpha_{n}\}_{n=0}^{\infty}$ associated with the or\-tho\-go\-nal polynomials with respect to $\mu$ are given by the relation
$$
\alpha_{n-1}=\overline{\tau}_{n-1}\left[\frac{1-2m_{n}-ic_{n}}{1-ic_{n}}\right], \quad n \geq 1,
$$
where $\tau_0 = 1$,  $\tau_{n}=\prod_{k=1}^{n}(1-ic_{k})/(1+ic_{k})$,  $n \geq 1$ and  $\{m_{n}\}_{n=0}^{\infty}$ is the minimal parameter sequence of $\{d_{n}\}_{n=1}^{\infty}$.
   In this manuscript we consider this relation and its consequences by imposing some restrictions of sign and periodicity  on the sequences $\{c_{n}\}_{n=1}^{\infty}$ and $\{m_{n}\}_{n=1}^{\infty}$. When the sequence $
   \{c_{n}\}_{n=1}^{\infty}$ is of alternating sign, we use information about the zeros of associated para-orthogonal polynomials to show that there is a gap in the support of the measure in the neighbourhood of $z= -1$.  Furthermore, we show that it is possible to ge\-nerate periodic Verblunsky coefficients by choosing periodic sequences $\{c_{n}\}_{n=1}^{\infty}$ and $\{m_{n}\}_{n=1}^{\infty}$  with the additional restriction  $c_{2n}=-c_{2n-1}, \, n\geq 1.$  We also give some results on periodic Verblunsky coefficients from the point of view of positive chain sequences. An example  is provided to illustrate the results obtained.
\end{abstract}

{\noindent}Keywords: Para-orthogonal polynomials, Probability measures, Periodic Verblunsky coefficients, Chain sequences, Alternating sign sequence. \\

{\noindent}2010 Mathematics Subject Classification: 42C05, 33C47.

\setcounter{equation}{0}
\section{Introduction}  \label{Sec-Intro}

Orthogonal polynomials on the unit circle (OPUC) have been commonly known as Szeg\H{o} polynomials in honor of Gábor Szeg\H{o} who introduced them in the first half of the 20th century. Because of their applications in quadrature rules, signal processing, operator and spectral theo\-ry and many other topics, these polynomials have received a lot of attention in recent years (see, for example, \cite{{BreuerSimon-2010},{Castil-Garza-Marcell-2011},{Costa-Ranga-Godoy-2011},{Kheif-Golins-Pehers-Yudit-2011},{Pehers-2011},{Pehers-Volb-Yudit-2011}, {Simanek-2012},{Tsujimoto-Zhedanov-2009}}). For many years a first hand text for an introduction to these polynomials has been the classical book \cite{Szego-book-1939} of Szeg\H{o}. However, for recent and more up to date texts on this subject we refer to the two volumes of Simon \cite{{Simon-book-p1},{Simon-book-p2}}. For further interesting reading on this subject we refer  to Chapter 8 of Ismail's recent book \cite{Ismail-book}.

Given a nontrivial probability measure $\mu(z)=\mu(e^{i\theta})$ on the unit circle $\mathbb{T}=\{z=e^{i\theta} : 0 \leq \theta \leq 2\pi \},$ the associated sequence of OPUC $\{\phi_n\}$ are those with the property
\[
\begin{array}l
\dsp \int_{\mathbb{T}} \bar{z}^j \phi_n(z) d\mu(z) = \int_{0}^{2\pi} e^{-ij\theta} \phi_n(e^{i\theta}) d\mu(e^{i\theta}) = 0, \quad 0 \leq j \leq n-1, \quad n \geq 1.
\end{array}
\]
\noindent  Letting  $\kappa_{n}^{-2} = \|\phi_n\|^2 = \int_{\mathbb{T}} |\phi_n(z)|^2 d\mu(z)$, the orthonormal polynomials on the unit circle are $\varphi_{n}(z) = \kappa_{n} \phi_n(z)$, $n \geq 0$.

The polynomials $\phi_n(z)$, $n \geq 0$, considered as monic polynomials,  satisfy the so called forward and backward recurrence relations, respectively,
\begin{equation} \label{Szego-A-RR}
\begin{array}l
  \phi_n(z) =  z \phi_{n-1}(z) - \overline{\alpha}_{n-1}\, \phi_{n-1}^{\ast}(z), \\[1.5ex]
  \phi_n(z) = (1 - |\alpha_{n-1}|^2) z \phi_{n-1}(z) - \overline{\alpha}_{n-1} \phi_n^{\ast}(z),
\end{array}
n \geq 1,
\end{equation}
where $\alpha_{n-1} = - \overline{\phi_n(0)}$ and $\phi_n^{\ast}(z) = z^{n} \overline{\phi_n(1/\bar{z})}$ denotes the reversed (reciprocal) polynomial of $\phi_n(z)$.   The numbers $\alpha_{n}$, in recent years, have been referred to  as Verblunsky coefficients.  It is known that these coefficients  are such that $|\alpha_n| < 1$, $n \geq 0$. Moreover, the OPUC and the associated measure are completely determined from these coefficients  (see for example \cite{Simon-book-p1}, Theorem 1.7.11).  A very nice and short constructive proof of this last statement can be found in  \cite{Nevai}.

It was shown in \cite{Costa-Felix-Ranga-2013}  that given any nontrivial probability measure on the unit circle, then corresponding to this measure there exists a pair of real sequences $\{c_{n}\}_{n=1}^{\infty}$ and $\{d_{n}\}_{n=1}^{\infty},$ where $\{d_{n}\}_{n=1}^{\infty}$ is also a positive chain sequence. In Theorem \ref{teo1} we have given the complete information regarding this statement and its reciprocal. To be precise, the sequences $\{c_{n}\}_{n=1}^{\infty}$ and $\{d_{n}\}_{n=1}^{\infty}$ are the coefficients of the three term recurrence formula 
\begin{eqnarray}\label{POPUC-A-RR}
R_{n+1}(z)=[(1+ic_{n+1})z+(1-ic_{n+1})]R_{n}(z)-4d_{n+1}zR_{n-1}(z), \ \ n\geq 1,
\end{eqnarray}
with $R_{0}(z)=1$ and $R_{1}(z)=(1+ic_{1})z+(1-ic_{1}),$ 
where  
\[
    R_n(z) = \frac{\prod_{j=1}^{n} \big[1- \tau_{j-1}\alpha_{j-1}\big]}{\prod_{j=1}^{n} \big[1-\mathcal{R}e(\tau_{j-1}\alpha_{j-1})\big]}\frac{z\phi_n(z) - \tau_{n} \phi_n^{\ast}(z)}{z-1},
\] 
with  $\tau_{n} = \phi_n(1)/\phi_n^{\ast}(1)$, $n \geq 0$. 

From the sequences $\{c_{n}\}_{n=1}^{\infty}$ and $\{d_{n}\}_{n=1}^{\infty},$ it is possible to recover the associated probability measure using certain rational functions that follow from the recurrence formula  (\ref{POPUC-A-RR}). In \cite{Castillo-Costa-Ranga-Veronese-2014}, using standard arguments involving continued fractions, series expansions at infinity and at the origin, and  Helly's Selection Theorem,  the associated measure $\mu$ is given as a limit of a subsequence of discrete  measures $\psi_{n}(e^{i\theta})$ whose pure points (those different from $z=1$)  are exactly the zeros of $R_{n}(z)$.  Results given in \cite{Castillo-Costa-Ranga-Veronese-2014} enable us to give information about the support of the measure $\mu$ by analysing  the zeros $z_{n,j} = e^{i\theta_{n,j}}$, $j = 1,2 \ldots, n$,  of $R_{n}(z),$ or, equivalently, by analysing the zeros of the functions $\mathcal{W}_{n}(x)$, given by 
\begin{equation} \label{Eq-for-Wn}
    \mathcal{W}_{n}(x) = 2^{-n} e^{-in\theta/2} R_{n}(e^{i\theta}), \quad n \geq 0,
\end{equation}
where $x = \cos(\theta/2)$. The sequence of functions $\{\mathcal{W}_{n}\}_{n=0}^{\infty}$ satisfy the three term recurrence formula (see \cite{BracMcCabPerezRanga-MCOM2015, DimRan-2013}) 
\begin{equation} \label{TTRR-for-Wn}
    \mathcal{W}_{n+1}(x) = \left(x - c_{n+1}\sqrt{1-x^2}\right)\mathcal{W}_{n}(x) - d_{n+1}\,\mathcal{W}_{n-1}(x), \quad n \geq 1,
\end{equation}
with $\mathcal{W}_{0}(x) = 1$ and $\mathcal{W}_{1}(x) = x - c_{1}\sqrt{1-x^2}$.  

For any $n \geq 1$, $\mathcal{W}_{n}(x)$ has exactly $n$ distinct zeros $x_{n,j} = \cos(\theta_{n,j}/2)$, $j = 1,2 \ldots, n$, in $(-1,1)$.  
We mention that the proof given in \cite{DimRan-2013} for the interlacing property of the zeros of $R_n(z)$ and $R_{n+1}(z)$ is by proving the interlacing property 
\begin{equation}  \label{Eq-InterlacingZeros-Wn}
    -1 < x_{n+1,n+1} < x_{n,n} < x_{n+1,n} < \cdots < x_{n,1} < x_{n+1,1} < 1, \quad n \geq 1,
\end{equation}
for the zeros of $\mathcal{W}_{n}(x)$ and $\mathcal{W}_{n+1}(x)$ using the three term recurrence formula (\ref{TTRR-for-Wn}).  

The aim of this manuscript is to study  sequences of Verblunsky coefficients where the related sequences  $\{c_{n}\}_{n=1}^{\infty}$ and $\{m_{n}\}_{n=1}^{\infty}$  have restrictions of sign and periodicity. We show that, under certain conditions, it is possible to estimate the support of the associated measure  and to get periodic Verblunsky coefficients. Furthermore, we discuss some geometric aspects related to these restrictions.

 This manuscript is organized as follows. In Section \ref{Sec-PrelimResults} we give a summary of all required theoretical results. Section \ref{alternatingmeasures} deals with the results concerning measures for which the associated sequence $\{c_{n}\}_{n=1}^{\infty}$ has the alternating sign property, namely, $c_{n}=(-1)^{n}\tilde{c}_{n}$, for $n \geq 1,$ where $\tilde{c}_{n}$ is a positive (or negative) sequence of real numbers. In Section \ref{periodicmeasures} relations with periodic Verblunsky coefficients are considered. Finally, in Section \ref{example} we give an example to illustrate the results obtained.

\setcounter{equation}{0}
\section{Some preliminary results}   \label{Sec-PrelimResults}

In this section we  present some results concerning nontrivial probability measures and positive chain sequences (for more details on chain sequences we refer to  \cite{Chihara-book} and  \cite{Wall}). Furthermore, some results about periodic Verblunsky coefficients are presented.

We begin with two theorems established in \cite{Costa-Felix-Ranga-2013}. The first theorem provides a characterization for nontrivial probability measures in terms of two sequences $\{c_{n}\}_{n=1}^{\infty}$ and $\{d_{n}\}_{n=1}^{\infty}.$ 

\begin{theo}\label{teo1}{\rm(a)} Given a nontrivial probability measure $\mu$ on the unit circle, then associated with it there exists an unique pair of real sequences $\left\{\{c_{n}\}_{n=1}^{\infty}, \{d_{n}\}_{n=1}^{\infty}\right\},$ where $\{d_{n}\}_{n=1}^{\infty}$ is also a positive chain sequence. Specifically, if $\{\alpha_{n}\}_{n=0}^{\infty}$ is the associated sequence of Verblunsky coefficients and if the sequence $\tau_{n}$ is such that
$$\tau_{0}=1 \quad \mbox{and} \quad \tau_{n}=\tau_{n-1}\frac{1-\overline{\tau}_{n-1}\overline{\alpha}_{n-1}}{1-\tau_{n-1}\alpha_{n-1}}, \quad n\geq1,$$
then $m_{0}=0,$ 
$$c_{n}=\frac{-Im(\tau_{n-1}\alpha_{n-1})}{1-Re(\tau_{n-1}\alpha_{n-1})} \quad \mbox{and} \quad m_{n}=\frac{1}{2}\frac{|1-\tau_{n-1}\alpha_{n-1}|^{2}}{[1-Re(\tau_{n-1}\alpha_{n-1})]}, \quad n\geq1,$$
where $\{m_{n}\}_{n=0}^{\infty}$ is the minimal parameter sequence of $\{d_{n}\}_{n=1}^{\infty}.$ Moreover, the maximal parameter sequence $\{M_{n}\}_{n=0}^{\infty}$ of $\{d_{n}\}_{n=1}^{\infty}$ is such that $M_{0}$ is the value of the jump in the measure at $z=1.$

\vspace{0.2cm}
\noindent{\rm(b)} Conversely, given a pair of real sequences $\left\{\{c_{n}\}_{n=1}^{\infty}, \{d_{n}\}_{n=1}^{\infty}\right\},$ where $\{d_{n}\}_{n=1}^{\infty}$ is also a positive chain sequence then associated with this pair there exists an unique nontrivial probability measure $\mu$ supported on the unit circle. Specifically, if $\{m_{n}\}_{n=0}^{\infty}$ is the minimal parameter sequence of $\{d_{n}\}_{n=1}^{\infty},$ then $\tau_{0}=1,$ 
\begin{equation}\label{verblunsky}
\tau_{n-1}\alpha_{n-1}=\frac{1-2m_{n}-ic_{n}}{1-ic_{n}} \quad \mbox{and} \quad \tau_{n}=\frac{1-ic_{n}}{1+ic_{n}}\tau_{n-1}, \quad n \geq 1.
\end{equation}
Moreover, the measure has a jump $M_{0}$ at $z=1,$ where $\{M_{n}\}_{n=0}^{\infty}$ is the maximal parameter sequence of $\{d_{n}\}_{n=1}^{\infty}.$
\end{theo}

The next  theorem gives information regarding the pure points. This  theorem is obtained as a consequence of Wall's criterion for maximal parameter sequence of positive chain sequences. 

\begin{theo}\label{teo2}
The probability measure $\mu$ has a pure point at $w\,(|w|=1)$ if, and only if,
 $$\sum_{n=1}^{\infty}\left[\prod_{j=1}^{n}\frac{|1-w\tau_{j-1}(w)\alpha_{j-1}|^{2}}{1-|\alpha_{j-1}|^{2}}\right]=\lambda(w)<\infty.$$
 Moreover, the size of the mass at the point $z=w$ is equal to $t=[1+\lambda(w)]^{-1}.$  Here, $\tau_{0}(w)=1$ and
\begin{equation}\label{taujw}
\tau_{j+1}(w)=\frac{\phi_{j+1}(w)}{\phi_{j+1}^{*}(w)}=\frac{w\tau_{j}(w)-\overline{\alpha}_{j}}{1-w\tau_{j}(w)\alpha_{j}}, \quad j\geq 0.
\end{equation}
\end{theo}

Now we discuss a result obtained in \cite{Castillo-Costa-Ranga-Veronese-2014}, which leads to a relation between the zeros of the polynomials $R_{n}(z)$ and the measure associated with the pair of sequences $\left\{\{c_{n}\}_{n=1}^{\infty}, \{d_{n}\}_{n=1}^{\infty}\right\}$.

Consider the new sequence of polynomials  $\{Q_{n}\}$ satisfying
$$Q_{n+1}(z)=[(1+ic_{n+1})z+(1-ic_{n+1})]Q_{n}(z)-4d_{n+1}zQ_{n-1}(z), \ \ n\geq 1,$$
with $Q_{0}(z)=0$ and $Q_{1}(z)=2d_{1}.$

Let $z_{n,j}=e^{i\theta_{n,j}}$ be the zeros of $R_{n}(z),$ $\lambda_{n,0}=1-\frac{Q_{n}(1)}{R_{n}(1)}$ and $\lambda_{n,j}=\frac{Q_{n}(z_{n,j})}{(1-z_{n,j})R_{n}^{\prime}(z_{n,j})},$ with
$j \in \{1,2,...,n\}.$ Thus, as shown in \cite{Castillo-Costa-Ranga-Veronese-2014},   $\sum_{j=0}^{n}\lambda_{n,j}=1$ and also $\lambda_{n,j}>0\, ,j=0,1,2...,n.$

In addition, if we define the sequence of step-functions $\psi_{n}(e^{i\theta}),n\geq1,$ on $[0,2\pi]$ by

\begin{equation}
\psi_{n}(e^{i\theta})=\left\{\begin{array}{ll}
                        0, & \theta=0, \\
                         \lambda_{n,0}, & 0<\theta\leq\theta_{n,1}, \\
                        \sum_{j=0}^{k}\lambda_{n,j}, & \theta_{n,k}<\theta\leq\theta_{n,k+1}, \quad k=1,2,...,n-1,\\
                        1, &  \theta_{n,n}<\theta\leq 2\pi
                      \end{array}\right.
\end{equation}
then by Helly Selection Theorem  a subsequence of $\psi_{n}(e^{i\theta})$ converges  to the measure $\mu(e^{i\theta})$ associated with the pair
$\left\{\{c_{n}\}_{n=1}^{\infty}, \{d_{n}\}_{n=1}^{\infty}\right\}$ as established in Theorem \ref{teo1}.

As an immediate consequence of this result, we can state the following.

\begin{theo}\label{suporte} Let $\left\{\{c_{n}\}_{n=1}^{\infty}, \{d_{n}\}_{n=1}^{\infty}\right\}$ be a pair of real sequences with $\{d_{n}\}_{n=1}^{\infty}$ a positive chain sequence. Moreover, let $R_{n}(z)$ be the sequence of polynomials given by {\rm(\ref{POPUC-A-RR})} and $\mu$ be the measure associated with  this pair of sequences. In addition, suppose that the zeros of $R_{n}(z)$ lie on an closed arc $\mathcal{B}$ of the unit circle, for $n\geq1.$ Then, the support of the measure $\mu$  lie within $\mathcal{B}\cup\{1\}.$
\end{theo}

 We now present a review of basic results on measures with periodic Verblunsky coefficients. For more details regarding these results we refer to   \cite{{Geronimus-1941}, {Pehers-Stein-1996}, {Pehers-Stein-1997}, {Simon-book-p2}}.

Let $\{\alpha_{n}\}_{n=0}^{\infty}$ be  a $p$-periodic sequence ($\alpha_{n+p}=\alpha_{n}, \, n\geq0$) of Verblunsky coefficients associated with the measure denoted by $\mu^{(p)}$ (here, $p$ is a fixed natural number).
 Consider the discriminant function $\Delta(z)=z^{-p/2}{\rm Tr}(T_{p}(z))$ where
  \begin{equation}\label{transferencia}
 T_{p}(z)=A(\alpha_{p-1},z)\ldots A(\alpha_{0},z),
 \end{equation}
 \begin{equation} \label{matrizA}
 A(\alpha_{j},z)= (1-|\alpha_{j}|^2)^{-1/2}\left(\begin{array}{cc}
                        z &  -\overline{\alpha}_{j}\\
                         -\alpha_{j} z & 1
           \end{array}\right), \quad j=0, \ldots, p-1,
 \end{equation}
 and ${\rm Tr}(T_{p}(z))$ denotes the trace of $T_{p}(z).$

 It is well known that all the $p$ distinct solutions of the equation $\Delta(z)=2$, which we denote by $z_{1}^{+}, \ldots , z_{p}^{+}$, lie on the unit circle $\mathbb{T}.$ In the same way, the $p$ distinct solutions of the equation $\Delta(z)=-2$, denoted by $z_{1}^{-}, \ldots z_{p}^{-}$, also lie on $\mathbb{T}.$ Using these solutions it is possible to show that the unit circle can be decomposed into $2p$ alternating sets $G_{1},B_{1},G_{2},\ldots,B_{p}$ with each gap, $G_{j},$ open and each band, $B_{j},$ closed. Moreover, each band $B_{j}$ is given by $B_{j} = \{z\in\mathbb{T}\, | \,\arg(z_{j}^{\sigma_{j}}) \leq \arg(z) \leq \arg(z_{j}^{-\sigma_{j}})\}$ with $\sigma_{j}=(-1)^{j+1},\, j=1, 2, \ldots, p.$

 Now we mention four fundamental results (see \cite[Chapter $11$]{Simon-book-p2}) which give a completely characterization of probability measures on the unit circle associated with periodic Verblunsky coefficients. 
 The  first result provides information about the absolutely continuous part and the singular part of the measure.

 \begin{theo}Let $\{\alpha_{j}\}_{j=0}^{\infty}$ be a sequence of Verblunsky coefficients of period $p$ and let $d\mu^{(p)}=w(\theta)\frac{d\theta}{2\pi}+d\mu^{(p)}_{s}$ be the associated probability measure. Then, if $B_{1}, \ldots, B_{p}$ are the corresponding bands we have that $\cup B_{j}$ is the essential support of the a.c. spectrum and $d\mu^{(p)}_{s}[\cup B_{j}]=\emptyset.$ Moreover, in each disjoint open arc on $\mathbb{T}\backslash \cup_{j=1}^{p} B_{j},$ $\mu^{(p)}$ has either no support or a single pure point.
  \end{theo}

  The next theorem provides information  about the associated weight function $w(\theta).$
  \begin{theo}\label{abscont}
  Let $\{\alpha_{j}\}_{j=0}^{\infty}$ be a sequence of Verblunsky coefficients of period $p$ and let $d\mu^{(p)}=w(\theta)\frac{d\theta}{2\pi}+d\mu^{(p)}_{s}$ be the associated measure. Then, for $e^{i\theta} \in \cup B_{j},$
\begin{equation}
w(\theta)=\frac{\sqrt{4-\Delta^{2}(e^{i\theta})}}{2|{\rm Im}(e^{-ip\theta/2})\varphi_{p}(e^{i\theta})|}\,.
\end{equation}
In particular,

\begin{enumerate}

\item [{\rm(i)}] On $\cup B_{j}^{int}, w(\theta)>0.$
\item [{\rm(ii)}]At an edge of a band that is by a closed gap (a gap which is empty), $w(\theta)>0.$
 \item [{\rm(iii)}] At an edge, $\theta_{0},$ of a band that is by an open gap, $w(\theta)\thicksim c(\theta-\theta_{0})^{\frac{1}{2}}$ if $\varphi_{p}^{*}(e^{i\theta_{0}})-\varphi_{p}(e^{i\theta_{0}})\neq 0.$
 \item [{\rm(iv)}]At an edge, $\theta_{0},$ of a band that is by an open gap, $w(\theta)\thicksim c(\theta-\theta_{0})^{-\frac{1}{2}}$ if $\varphi_{p}^{*}(e^{i\theta_{0}})-\varphi_{p}(e^{i\theta_{0}})= 0.$
\end{enumerate}
  \end{theo}

  Finally, the following two theorems lead to a complete characterization for the pure points of the measure.

  \begin{theo}\label{puro1}
  Let $\{\alpha_{j}\}_{j=0}^{\infty}$ be a sequence of Verblunsky coefficients of period $p.$ Then
  \begin{equation}
  \pi(z)=\varphi_{p}^{*}(z)-\varphi_{p}(z)
\end{equation}
 has all its zeros in the set of gap closures, one in each gap closure.
  \end{theo}

  \begin{theo}\label{puro2}
  Let $\{\alpha_{j}\}_{j=0}^{\infty}$ be a sequence of Verblunsky coefficients of period $p$ and let $\mu^{(p)}$ the associated measure. Let $\theta_{0}$ be a point in a gap closure where $\varphi_{p}^{*}(e^{i\theta_{0}})-\varphi_{p}(e^{i\theta_{0}})=0.$ Then, either $\mu^{(p)}$ has no pure point in the gap or else it has a pure point at $z=e^{i\theta_{0}}.$
   \end{theo}

   The  results on periodic Verblunsky coefficients presented above will be used in sections \ref{periodicmeasures} and \ref{example}.

\setcounter{equation}{0}
\section{On measures associated with alternating sign sequences $\{c_{n}\}$ }   \label{alternatingmeasures}

First, we provide two lemmas that will be useful to derive the subsequent  results.

\begin{lema}\label{equivalencia}
Let $\mathcal{W}_{n}(x)$ satisfying {\rm(\ref{TTRR-for-Wn})} and  $R_n(z)$ satisfying {\rm(\ref{POPUC-A-RR})}. Then, the following statements are equivalent:
\begin{enumerate}
\item[{\rm(i)}] $c_n=(-1)^nc$, $n\geq 1$  and $c\in \mathbb{R}$;
\item[{\rm(ii)}] For $n\geq0$, $R_{2n}(z)$ has real coefficients and {\small$R_{2n+1}(z)=\left[(1-ic)z+(1+ic)\right]\widetilde{R}_{2n}(z)$}, where $\widetilde{R}_{0}(x)=1$ and $\widetilde{R}_{2n}(z)$ is also a polynomial with real coefficients;
\item[{\rm(iii)}] For $n\geq0$, $\mathcal{W}_{2n}(x)$ is an even polynomial of degree $2n$ and {\small$\mathcal{W}_{2n+1}(x)=\left(x+c\sqrt{1-x^2}\right)\widetilde{\mathcal{W}}_{2n}(x)$} with  $\widetilde{\mathcal{W}}_{0}(x)=1$ and $\widetilde{\mathcal{W}}_{2n}(x)$ an even polynomial of degree $2n$. 
\end{enumerate}
\end{lema}

{\proof}
(i)$\Rightarrow$(ii) Since $R_{0}(z)=1$ and    $R_{1}(z)=\left[(1-ic)z+(1+ic)\right]\widetilde{R}_{0}(z),$
 with $\widetilde{R}_{0}(z)=1,$ it follows that the result holds for $n=0.$ Furthermore, if (ii) holds for $n=k\in\mathbb{N}$ then, from the three term recurrence relation (\ref{POPUC-A-RR}), we obtain
\begin{eqnarray*}
R_{2(k+1)}(z)&=&\left[(1+ic)z+(1-ic)\right]\left[(1-ic)z+(1+ic)\right]\widetilde{R}_{2k}(z)-4d_{2k+2}zR_{2k}(z)\\
&=& \left[(1+c^2)z^2+2(1-c^2)z+(1+c^2)\right]\widetilde{R}_{2k}(z)-4d_{2k+2}zR_{2k}(z).
\end{eqnarray*}
Consequently, since we are assuming that $\widetilde{R}_{2k}(z)$ and $R_{2k}(z)$ are polynomials with real coefficients, we conclude that $R_{2(k+1)}(z)$ also has real coefficients. Moreover, using again (\ref{POPUC-A-RR}), we can see that
\begin{eqnarray*}
R_{2(k+1)+1}(z)&=&\left[(1-ic)z+(1+ic)\right]R_{2(k+1)}(z)-4d_{2k+3}z\left[(1-ic)z+(1+ic)\right]\widetilde{R}_{2k}(z)\\
&=&\left[(1-ic)z+(1+ic)\right]\widetilde{R}_{2(k+1)}(z),
\end{eqnarray*}
where $\widetilde{R}_{2(k+1)}(z)=R_{2(k+1)}(z)-4d_{2k+3}z\widetilde{R}_{2k}(z)$ is also a polynomial with real coefficients, once $R_{2(k+1)}(z)$ and $\widetilde{R}_{2k}(z)$ have real coefficients. Therefore, using mathematical induction, we conclude that the statement (ii) holds for all $n\geq0$.

  (ii)$\Rightarrow$(iii) By (ii), $R_{2n}(z)$ has real coefficients for $n\geq 0$. Moreover, from the three term recurrence relation (\ref{POPUC-A-RR}) we have that $R_{2n}(z)$ is a self-inversive polynomial, i.e., $R_{2n}^{\ast}(z) = z^{2n} \overline{R_{2n}(1/\bar{z})}=R_{2n}(z)$. Therefore, it follows (see \cite[Lemma $2.1$]{BracMcCabPerezRanga-MCOM2015}) that $W_{2n}(x)$ is an even polynomial of degree $2n$ in the variable $x=\cos \theta/2$. Similarly, since $\widetilde{R}_{2n}(z)$ is  also a self-inversive polynomial with real coefficients, we have that $\widetilde{\mathcal{W}}_{2n}(x)=(4e^{i\theta})^{-n/2} \widetilde{R}_{2n}(e^{i\theta})$ is an even polynomial of degree $2n$. Then, since $R_{2n+1}(z)=\left[(1-ic)z+(1+ic)\right]\widetilde{R}_{2n}(z),$ from the relation (\ref{Eq-for-Wn}) it follows that $\mathcal{W}_{2n+1}(x)=(x+c\sqrt{1-x^2})\widetilde{\mathcal{W}}_{2n}(x).$  
 
  (iii)$\Rightarrow$(i) Using the assumption (iii) and the three term recurrence relation (\ref{TTRR-for-Wn}), we have, for $s\geq1$,
\begin{equation}\label{w1j}
\mathcal{W}_{2s}(x)=(x - c_{2s}\sqrt{1-x^2})(x+c\sqrt{1-x^2})\widetilde{\mathcal{W}}_{2s-2}(x) - d_{2s}\,\mathcal{W}_{2s-2}(x)
\end{equation}
and
\begin{equation}\label{w2j}
(x+c\sqrt{1-x^2})\widetilde{\mathcal{W}}_{2s}(x)=(x - c_{2s+1}\sqrt{1-x^2})\mathcal{W}_{2s}(x) - d_{2s+1}(x+c\sqrt{1-x^2})\widetilde{\mathcal{W}}_{2s-2}(x).
\end{equation}
Hence, since $\mathcal{W}_{2s}(x)$, $\mathcal{W}_{2s-2}(x)$, $\widetilde{\mathcal{W}}_{2s}(x)$ and $\widetilde{\mathcal{W}}_{2s-2}(x)$ are even polynomials, we can use the relations (\ref{w1j}) and (\ref{w2j}) to conclude that $c_{2s}=c$ and $c_{2s+1}=-c$, for $s\geq1$. Moreover, using the assumption (iii) and the definition of $\mathcal{W}_{1}(x)$ it is easy to see that $c_{1}=-c.$
\fim

\vspace{0.3cm}
  Consider now the polynomials $\hat{R}_{n}(z)$ satisfying
  \begin{equation}\label{rchapeu}
  \hat{R}_{n+1}(z)=[(1+i\hat{c}_{n+1})z+(1-i\hat{c}_{n+1})]\hat{R}_{n}(z)-4d_{n+1}z\hat{R}_{n-1}(z), \ \ n\geq 1,\end{equation}
      with $\hat{R}_{0}(z)=1,$ $\hat{R}_{1}(z)=(1+i\hat{c}_{1})z+(1-i\hat{c}_{1})$ and $\hat{c}_{n}=-c_{n}.$

      The following lemma gives the relation between the polynomials $R_{n}(z)$ and $\hat{R}_{n}(z).$
        \begin{lema}\label{lemarrchapeu}
Let $R_{n}(z)$ satisfying {\rm(\ref{POPUC-A-RR})} and $\hat{R}_{n}(z)$ satisfying {\rm(\ref{rchapeu})}. Then, the following holds
$$R_{n}(z)=\overline{\hat{R}_{n}(\bar{z})}, \quad n=0, 1, 2, \ldots\, \ .$$
  \end{lema}

  {\proof}
  The proof can be given by mathematical induction. Clearly, the result holds for $n=0$ and $n=1.$ Suppose that the result holds for $n=0, 1, \ldots, k.$ Then, from the recurrence relations {\rm(\ref{POPUC-A-RR})} and {\rm(\ref{rchapeu})}, we have
  \begin{eqnarray*}
  \hat{R}_{k+1}(\bar{z})&=& [(1-ic_{k+1})\bar{z}+(1+ic_{k+1})]\hat{R}_{k}(\bar{z})-4d_{k+1}\bar{z}\hat{R}_{k-1}(\bar{z})\\
  &=& [(1-ic_{k+1})\bar{z}+(1+ic_{k+1})]\overline{{R}_{k}(z)}-4d_{k+1}\bar{z}\overline{{R}_{k-1}(z).}
  \end{eqnarray*}
  Hence, the result follows by taking the complex conjugate on the above relation.
  \fim

  \vspace{0.3cm}
  Observe that the Lemma \ref{lemarrchapeu} provides also a relation between the zeros of the polynomials $R_{n}(z)$ and the zeros of $\hat{R}_{n}(z),$ namely, if $z_{n,j}$ is a zero of $R_{n}(z)$ then $\overline{z}_{n,j}$ is a zero of $\hat{R}_{n}(z).$

  Now we consider the problem of giving estimates for the support of measures whose sequences $\{c_{n}\}_{n=1}^{\infty}$ are of alternating sign. We start with the case $c_n=(-1)^nc$, where $c\in \mathbb{R}.$
  
  Let $\mathcal{C}_{1}$ and $\mathcal{C}_{2}$ be closed arcs on the unit circle given by
 $$\mathcal{C}_{1}=\left\{z\in\mathbb{T}\, :\, 0\leq \arg(z) \leq \theta_c\right\} 
 \quad \mbox{and} \quad \mathcal{C}_{2}=\left\{z\in\mathbb{T}\, :\, 2\pi - \theta_c \leq \arg(z) \leq 2\pi\right\},$$
where $\theta_c=\arccos \left(\frac{c^{2}-1}{c^{2}+1}\right)\in[0, \pi].$
Then, we can state the following.
  \begin{theo}\label{teoc}
  Let $\mu$ be the probability measure on the unit circle associated with the pair of sequences $\left\{\{c_{n}\}_{n=1}^{\infty}, \{d_{n}\}_{n=1}^{\infty}\right\}$ where $c_n=(-1)^nc,\, {\rm c \in \mathbb{R}}$ and $\{d_{n}\}_{n=1}^{\infty}$ is a positive chain sequence. Then, the support of $\mu$ lie within $\mathcal{C}_{1} \cup \mathcal{C}_{2}.$
  \end{theo}

  {\proof}
  Without loss of generality we assume that
   $c \geq 0.$
     Consider the polynomials $R_{n}(z)$ given by (\ref{POPUC-A-RR}). If we show that all zeros of $R_{n}(z)$ lie on $\mathcal{C}_{1} \cup \mathcal{C}_{2},$ then from Theorem \ref{suporte} we obtain the desired result.
  To show this, we use the functions $\mathcal{W}_{n}(x)$ defined in (\ref{TTRR-for-Wn}) which are associated to the polynomials $R_{n}(z).$

  By Lemma \ref{equivalencia} we have that $\mathcal{W}_{2n+1}(x)=(x+c\sqrt{1-x^2})\widetilde{\mathcal{W}}_{2n}(x)$ with  $\widetilde{\mathcal{W}}_{2n}(x)$ an even polynomial of degree $2n.$ Moreover, $\mathcal{W}_{2n}(x)$ is also an even polynomial of degree $2n.$ This means that $\frac{-c}{\sqrt{1+c^{2}}}$ is always a zero of $\mathcal{W}_{2n+1}(x)$ and the other $2n$ zeros of these functions have a symmetry about the origin. Likewise, $\mathcal{W}_{2n}(x)$ being an even polynomial,  all of their zeros are symmetric with respect to the origin.

  Therefore, from the symmetry of the zeros observed above and taking into account the interlacing property (\ref{Eq-InterlacingZeros-Wn}) it follows that all zeros of $\mathcal{W}_{n}(x)$ lie in $\left(-1,\frac{-c}{\sqrt{1+c^{2}}} \right] \cup \left[\frac{c}{\sqrt{1+c^{2}}},1\right).$

  Finally, if we denote the zeros of $\mathcal{W}_{n}(x)$ by $x_{n,j}$ and the zeros of $R_{n}(z)$ by $z_{n,j},$ then they are related by $x_{n,j}=\cos\left(\frac{\theta_{n,j}}{2}\right)$ where $z_{n,j}=e^{i\theta_{n,j}}, \, j=1, 2, \ldots, n.$ This shows that $R_{n}(z)$ has all of its zeros on $\mathcal{C}_{1} \cup \mathcal{C}_{2}.$
  \fim\vspace{0.3cm}

 Notice that Theorem \ref{teoc} leads to an estimative for the support of the measure in the case where $c_{n}=(-1)^{n}\tilde{c}_{n}$ and $\tilde{c}_{n}$ is a constant sequence. We use this initial estimative as motivation to obtain a more general result.

  \begin{theo}\label{teo-geral-alternantes}
  Let $\mu$ be the probability measure on the unit circle associated with the pair of sequences $\left\{\{c_{n}\}_{n=1}^{\infty}, \{d_{n}\}_{n=1}^{\infty}\right\}$ where $c_n=(-1)^{n}\tilde{c}_{n},\, \tilde{c}_{n}\geq c>0$ and $\{d_{n}\}_{n=1}^{\infty}$ is a positive chain sequence. Then, the support of $\mu$ lie within $\mathcal{C}_{1} \cup \mathcal{C}_{2}.$
  \end{theo}

  {\proof}
  Firstly, notice that for $\varepsilon$ with $0<\varepsilon<c,$  we have $\tilde{c}_{n}\geq c > c_{\varepsilon}>0,$ where $c_{\varepsilon}=c-\varepsilon$.
   Since $\tilde{c}_{n}>c_{\varepsilon}>0,$ for $x_{0}=\frac{-c_{\varepsilon}}{\sqrt{1+c_{\varepsilon}^{2}}}$ and $x_{1}=\frac{c_{\varepsilon}}{\sqrt{1+c_{\varepsilon}^{2}}}$, one can observe that, for $n\geq 1$,
   \begin{equation}\label{sign1}
      \mbox{sign}\left(x_{0}-\tilde{c}_{n}\sqrt{1-x_{0}^2}\right)=\mbox{sign}
      \left(x_{1}-\tilde{c}_{n}\sqrt{1-x_{1}^2}\right)=-1
   \end{equation}
      and
   \begin{equation}\label{sign2}
      \mbox{sign}\left(x_{0}+\tilde{c}_{n}\sqrt{1-x_{0}^2}\right)=\mbox{sign}
      \left(x_{1}+\tilde{c}_{n}\sqrt{1-x_{1}^2}\right)=1.
   \end{equation}
  
   Now, we will show that all zeros of $\mathcal{W}_{n}(x)$ lie in $\left(-1,x_{0} \right] \cup \left[x_{1},1\right)$. Noticing that $\frac{-\tilde{c}_{1}}{\sqrt{1+\tilde{c}_{1}^{2}}}$ is the only zero of $\mathcal{W}_{1}(x)$ in $(-1,1)$ and that $\tilde{c}_{1}>c_{\varepsilon}$, it follows that the result is valid for $n=1.$ In addition, from (\ref{sign2}), $\mbox{sign}(\mathcal{W}_{1}(x_{j}))=1, \,j\in \{0,1\}.$

   Hence, from the three term recurrence relation (\ref{TTRR-for-Wn}) for $\mathcal{W}_{2}(x)$ and from (\ref{sign1}) we conclude that $\mbox{sign}(\mathcal{W}_{2}(x_{j}))=-1, \,j\in \{0,1\}.$ Suppose that there exists at least one zero of $\mathcal{W}_{2}(x)$ inside the interval $(x_{0},x_{1}).$ Then, since $\mathcal{W}_{2}(x_{0})<0$ and $\mathcal{W}_{2}(x_{1})<0$ we conclude that $\mathcal{W}_{2}(x)$ has two zeros in $(x_{0},x_{1}).$ But this cannot happen because the only zero  of $\mathcal{W}_{1}(x)$ is outside of $(x_{0},x_{1})$ and the two zeros of $\mathcal{W}_{2}(x)$ interlace with the zero of $\mathcal{W}_{1}(x)$. Thus, the result also holds for $n=2.$

   Again, from the recurrence relation for $\mathcal{W}_{3}(x),$ $\mbox{sign}(\mathcal{W}_{2}(x_{j}))= -1, \, \mbox{sign}(\mathcal{W}_{1}(x_{j}))=1, \, j\in \{0,1\},$ and by (\ref{sign2}) it follows that $\mbox{sign}(\mathcal{W}_{3}(x_{j}))=-1, \,j\in \{0,1\}.$ Hence, using the interlacing property for the zeros of $\mathcal{W}_{3}(x)$ and $\mathcal{W}_{2}(x),$ and the fact that there exist no zeros of $\mathcal{W}_{2}(x)$ in $(x_{0},x_{1}),$  it follows that $\mathcal{W}_{3}(x)$ cannot vanish in $(x_{0},x_{1}).$

   Continuing this procedure, by mathematical induction, one can easily see that
   $$\mbox{sign}(\mathcal{W}_{n}(x_{j}))=(-1)^{\left \lfloor n/2 \right \rfloor}, \ \,j\in \{0,1\}, \ \, n=0, 1, 2, \ldots$$
   and, by the same arguments used before, $\mathcal{W}_{n}(x)$, $n\geq1$, cannot vanish in $(x_{0},x_{1}).$ 
   
   Finally, by letting $\varepsilon\rightarrow0,$ we see that $\mathcal{W}_{n}(x)$ has all its zeros in $\left(-1,\frac{-c}{\sqrt{1+c^{2}}} \right] \cup \left[\frac{c}{\sqrt{1+c^{2}}},1\right)$ or, equivalently, $R_{n}(z), \, n \geq 1,$ has all of its zeros on $\mathcal{C}_{1} \cup \mathcal{C}_{2}.$ Now, the result follows by Theorem \ref{suporte}.
\fim

\begin{coro}
 Let $\mu$ be the probability measure on the unit circle associated with the pair of sequences $\left\{\{c_{n}\}_{n=1}^{\infty}, \{d_{n}\}_{n=1}^{\infty}\right\}$, where $c_n=(-1)^{n}\tilde{c}_{n},\, \tilde{c}_{n}\leq c<0$ and $\{d_{n}\}_{n=1}^{\infty}$ is a positive chain sequence. Then, the support of $\mu$ lie on $\mathcal{C}_{1} \cup \mathcal{C}_{2}.$
\end{coro}

{\proof}
First, one can observe that $-c_n=(-1)^{n}\left(-\tilde{c}_{n}\right),$ with $-\tilde{c}_{n}\geq -c>0.$ Hence, if $\hat{\mu}$ is the probability measure associated to the pair $\left\{\{-{c}_{n}\}_{n=1}^{\infty}, \{d_{n}\}_{n=1}^{\infty}\right\},$ from Theorem \ref{teo-geral-alternantes} it follows that $\hat{R}_{n}(z)$ given by (\ref{rchapeu}) has all zeros on $\mathcal{C}_{1} \cup \mathcal{C}_{2}$ and that $\hat{\mu}$ has its support within $\mathcal{C}_{1} \cup \mathcal{C}_{2}.$ Now the result is an immediate consequence of Lemma \ref{lemarrchapeu}.
\fim

\vspace{0.3cm}
Now we consider the measure $\mu$ associated the the pair $\left\{\{c_{n}\}_{n=1}^{\infty}, \{d_{n}\}_{n=1}^{\infty}\right\}$, where $\{c_{n}\}_{n=1}^{\infty}$  satisfy the condition $c_{2n}=-c_{2n-1}, \, n\geq 1.$

Starting from $\mu$ we desire to get a new measure $\tilde{\mu}$ associated with the pair $\left\{\{\tilde{c}_{n}\}_{n=1}^{\infty}, \{\tilde{d}_{n}\}_{n=1}^{\infty}\right\}$, where the sequence $\{\tilde{c}_{n}\}_{n=1}^{\infty}$ must satisfy the condition $\tilde{c}_{2n}=\tilde{c}_{2n-1}=c_{2n},\, n\geq 1.$

Let us consider the sequence of complex numbers $\{\beta_{n}\}_{n=1}^{\infty}$ given by
\begin{equation}\label{wn2}
\beta_{n}=-\left(\frac{1+ic_{2n}}{1-ic_{2n}}\right), \ \, n=1, 2, \ldots\, \ .
\end{equation}
\vspace{0.2cm}
The next theorem shows how to get the required measure $\tilde{\mu}$ from the measure $\mu.$
     \begin{theo} \label{teoprod}
      Let $\mu$ be the probability measure on the unit circle associated with the pair of sequences $\left\{\{c_{n}\}_{n=1}^{\infty}, \{d_{n}\}_{n=1}^{\infty}\right\}$ where  $c_{2n}=-c_{2n-1}, \, n \geq 1.$  Let $\{\beta_{n}\}_{n=1}^{\infty}$ be the sequence of complex numbers defined by {\rm(\ref{wn2})}. In addition, let  $\tilde{\mu}$ be the measure associated with the sequence of Verblunsky coefficients $\{\tilde{\alpha}_{n}\}_{n=0}^{\infty}$ given by
              \begin{eqnarray} \label{prod1}
              \tilde{\alpha}_{2n+1}=\left(\prod_{j=1}^{n+1}{\beta_{j}}^{2}\right)\alpha_{2n+1}
              \quad
       \mbox{and} \quad
              \tilde{\alpha}_{2n}=\left(\prod_{j=1}^{n}{\beta_{j}}^{2}\right)\beta_{n+1}\alpha_{2n}, \quad n=0, 1, 2, \ldots\,,
       \end{eqnarray}
       where $\{\alpha_{n}\}_{n=0}^{\infty}$ is the sequence of Verblunsky coefficients corresponding to $\mu.$
      If {\small $\left\{\{\tilde{c}_{n}\}_{n=1}^{\infty}, \{\tilde{d}_{n}\}_{n=1}^{\infty}\right\}$} is the pair of sequences associated with the measure $\tilde{\mu}$ and if $\{\tilde{m}_{n}\}_{n=0}^{\infty}$ is the minimal parameter sequence for $\{\tilde{d}_{n}\}_{n=1}^{\infty},$ then the following holds
$$ \tilde{c}_{2n}=\tilde{c}_{2n-1}=c_{2n}, \, \  \,\tilde{m}_{2n-1}=1-m_{2n-1} \quad \mbox{and} \quad \tilde{m}_{2n}=m_{2n}, \quad n=1, 2, \ldots\, \ .$$
            \end{theo}

            {\proof}  Using the assumption $c_{2n}=-c_{2n-1},$ we obtain
 \begin{eqnarray}\label{tn}
           \tau_{2n}=1 \quad \mbox{and} \quad \tau_{2n+1}=\frac{1+ic_{2n+2}}{1-ic_{2n+2}}, \quad n=0, 1, \ldots\, \ .
 \end{eqnarray}
           
            Hence, from (\ref{verblunsky}) and (\ref{tn}), we have
 \begin{eqnarray}\label{alphan}
           \alpha_{2n}=\frac{1-2m_{2n+1}+ic_{2n+2}}{1+ic_{2n+2}} \quad \mbox{and} \quad \alpha_{2n+1}=\frac{1-2m_{2n+2}-ic_{2n+2}}{1+ic_{2n+2}}, \quad n=0, 1, \ldots\, \ .
 \end{eqnarray}

Now let $\{\hat{m}_{n}\}_{n=0}^{\infty}$ be the minimal parameter sequence for a positive chain sequence  $\{\hat{d}_{n}\}_{n=1}^{\infty}$ and $\{\hat{\alpha}_{n}\}_{n=0}^{\infty}$ the Verblunsky coefficients of a probability measure on the unit circle, $\hat{\mu}$, associated with the pair of real sequences $\left\{\{\hat{c}_{n}\}_{n=1}^{\infty}, \{\hat{d}_{n}\}_{n=1}^{\infty}\right\}$, where
 \begin{eqnarray}\label{cn1}
\hat{c}_{2n}=\hat{c}_{2n-1}=c_{2n}, \, \  \,\hat{m}_{2n-1}=1-m_{2n-1} \quad \mbox{and} \quad \hat{m}_{2n}=m_{2n}, \quad n=1, 2, \ldots\, \ .
\end{eqnarray}

Using the relations (\ref{verblunsky}), (\ref{wn2}), (\ref{prod1}), (\ref{alphan}) and (\ref{cn1}), one can see that for $n=0, 1, \ldots\, \ ,$
 \begin{eqnarray*}
\hat{\alpha}_{2n+1}&=&\left(\frac{1+i\hat{c}_{2n+1}}{1-i\hat{c}_{2n+1}}\right)\left(\prod_{k=1}^{2n}\frac{1+i\hat{c}_{k}}{1-i\hat{c}_{k}}\right)\left[\frac{1-2\hat{m}_{2n+2}-i\hat{c}_{2n+2}}{1-i\hat{c}_{2n+2}}\right]\\
&=&\left(\frac{1+ic_{2n+2}}{1-ic_{2n+2}}\right)^2\left[\prod_{j=1}^{n}\left(\frac{1+ic_{2j}}{1-ic_{2j}}\right)^2\right]\left[\frac{1-2m_{2n+2}-ic_{2n+2}}{1+ic_{2n+2}}\right]=\left(\prod_{j=1}^{n+1}{\beta_{j}}^{2}\right)\alpha_{2n+1}\\
&=&\tilde{\alpha}_{2n+1}.
 \end{eqnarray*}
 
 Similarly, using again (\ref{verblunsky}), (\ref{wn2}), (\ref{prod1}), (\ref{alphan}) and (\ref{cn1}), we obtain for $n=0, 1, \ldots\, \ ,$
  \begin{eqnarray*}
\hat{\alpha}_{2n}&=&\left(\prod_{k=1}^{2n}\frac{1+i\hat{c}_{k}}{1-i\hat{c}_{k}}\right)\left[\frac{1-2\hat{m}_{2n+1}-i\hat{c}_{2n+1}}{1-i\hat{c}_{2n+1}}\right]\\
&=&\left[\prod_{j=1}^{n}\left(\frac{1+ic_{2j}}{1-ic_{2j}}\right)^2\right]\left[-\left(\frac{1+ic_{2n+2}}{1-ic_{2n+2}}\right)\right]\left[\frac{1-2m_{2n+1}+ic_{2n+2}}{1+ic_{2n+2}}\right]=\left(\prod_{j=1}^{n}{\beta_{j}}^{2}\right)\beta_{n+1}\alpha_{2n}\\
&=&\tilde{\alpha}_{2n}.
 \end{eqnarray*}

 Thus, $\tilde{\alpha}_{n}=\hat{\alpha}_{n}$ for $n\geq0$ and, consequently, $\tilde{\mu}=\hat{\mu}$. Hence, from the uniqueness of the pair  $\left\{\{\tilde{c}_{n}\}_{n=1}^{\infty}, \{\tilde{d}_{n}\}_{n=1}^{\infty}\right\}$ given by Theorem \ref{teo1}, we have $\tilde{m}_{0}=\hat{m}_{0}=0,$
 $$\tilde{c}_{n}=\hat{c}_{n}\quad\mbox{ and }\quad\tilde{m}_{n}=\hat{m}_{n}, \quad\quad n=1, 2, \ldots\,,$$
 which completes the proof of the theorem.
       \fim

       \begin{coro}\label{rotacao}
       Let $\mu$ be the probability measure on the unit circle associated with the pair of sequences $\left\{\{c_{n}\}_{n=1}^{\infty}, \{d_{n}\}_{n=1}^{\infty}\right\}$, where  $c_{n}=(-1)^{n}c, \, n \geq 1$ and $c \in \mathbb{R}. $ In addition, let $\beta=-\left(\frac{1+ic}{1-ic}\right)$ and $\tilde{\mu}(z)=\mu(\beta z)$ the measure associated with the pair $\left\{\{\tilde{c}_{n}\}_{n=1}^{\infty}, \{\tilde{d}_{n}\}_{n=1}^{\infty}\right\}.$ Then, for $ n\geq 1, \, \tilde{c}_{n}=c.$
       \end{coro}

       {\proof}
       First, notice that if $\tilde{\mu}(z)=\mu(\beta z),$ the corresponding Verblunsky coefficients are related by $\tilde{\alpha}_{n}=\beta^{n+1}\alpha_{n}, \, n \geq 0$ (see, for example, \cite{Costa-Felix-Ranga-2013}). Hence, the result follows  from Theorem \ref{teoprod} with $c_{2n}=c, \, n \geq 1.$
       \fim

\setcounter{equation}{0}
\section{Measures with periodic Verblunsky coefficients}   \label{periodicmeasures}

The first theorem in this section gives a characterization of measures with periodic Verblunsky coefficients in terms of the pair of real sequences  $\left\{\{c_{n}\}_{n=1}^{\infty}, \{d_{n}\}_{n=1}^{\infty}\right\}$, where $\{d_{n}\}_{n=1}^{\infty}$ is a positive chain sequence. Throughout in this section $b_{n}=1-2m_{n}, \, n \geq 1,$ where $\{m_{n}\}_{n=0}^{\infty}$ is the minimal parameter sequence of $\{d_{n}\}_{n=1}^{\infty}.
$

\begin{theo} \label{caracterizacao}
 Let $\mu$ be the probability measure on the unit circle associated with the pair of sequences $\left\{\{c_{n}\}_{n=1}^{\infty}, \{d_{n}\}_{n=1}^{\infty}\right\}.$ Then, the measure $\mu$ has periodic Verblunsky coefficients $\{\alpha_{n}\}_{n=0}^{\infty}$ of period $p$ if, and only if, for $n \geq 0,$
 \begin{equation} \label{argumento}
 \sum_{j=n+1}^{n+p}\arg\left(\frac{1+ic_{j}}{1-ic_{j}}\right)= \arg\left(\frac{b_{n+1}-ic_{n+1}}{1-ic_{n+1}}\right)-\arg\left(\frac{b_{n+p+1}-ic_{n+p+1}}{1-ic_{n+p+1}}\right)+2k_{n}\pi, \, \,  k_{n} \in \mathbb{Z}
  \end{equation}
  and
  \begin{equation}\label{modulo}
  \frac{b_{n+1}^{2}+c_{n+1}^{2}}{1+c_{n+1}^{2}}=\frac{b_{n+p+1}^{2}+c_{n+p+1}^{2}}{1+c_{n+p+1}^{2}}.
  \end{equation}
\end{theo}

{\proof}
First one can observe, from (\ref{verblunsky}), that for $n \geq 0$
\begin{eqnarray*}\alpha_{n+p}=\alpha_{n} \, \, &\Leftrightarrow& \,\, \overline{\tau}_{n+p}\left[ \frac{b_{n+p+1}-ic_{n+p+1}}{1-ic_{n+p+1}}\right]=\overline{\tau}_{n} \left[\frac{b_{n+1}-ic_{n+1}}{1-ic_{n+1}}\right] \,\, \\
&\Leftrightarrow& \,\, \left(\prod_{j=n+1}^{n+p}\frac{1+ic_{j}}{1-ic_{j}}\right) \left[\frac{b_{n+p+1}-ic_{n+p+1}}{1-ic_{n+p+1}}\right]= \left[\frac{b_{n+1}-ic_{n+1}}{1-ic_{n+1}}\right].
\end{eqnarray*}

 Now the result follows by comparing, respectively, the modulus and the argument of the numbers
 $$\displaystyle{\left(\prod_{j=n+1}^{n+p}\frac{1+ic_{j}}{1-ic_{j}}\right) \left[\frac{b_{n+p+1}-ic_{n+p+1}}{1-ic_{n+p+1}}\right]}\quad \mbox{and} \quad  \displaystyle{\left[\frac{b_{n+1}-ic_{n+1}}{1-ic_{n+1}}\right]}, \, \ n\geq 0.$$
 \fim\vspace{0.2cm}
 
 We say that $\mu$ is a symmetric measure if $d\mu(z)=-d\mu(1/z), \, z \in \mathbb{T}.$ From results established in \cite{Castillo-Costa-Ranga-Veronese-2014} one can observe that $\mu$ is symmetric if and only if $c_{n}=0$, $n \geq 1$, with $\{c_{n}\}_{n=1}^{\infty}$ given as in Theorem \ref{teo1}. Thus, as a consequence of Theorem \ref{caracterizacao}, we have the following result.

\begin{coro}\label{periodicaalternante}
Let $\mu$ be the probability measure on the unit circle associated with the pair of sequences $\left\{\{c_{n}\}_{n=1}^{\infty}, \{d_{n}\}_{n=1}^{\infty}\right\},$ where $\{c_{n}\}_{n=1}^{\infty}$ and $\{m_{n}\}_{n=1}^{\infty}$ are periodic sequences of period $p.$ In addition, suppose that $c_{2n}=-c_{2n-1}, \, n \geq 1.$ Then,
\begin{enumerate}
\item[{\rm(i)}] if $p$ is even the measure $\mu$ has $p\,-$periodic sequence of Verblunsky coefficients;
\item[{\rm(ii)}] if $p$ is odd, the measure $\mu$ is symmetric and has $p\,-$periodic sequence of Verblunsky coefficients.
\end{enumerate}
\end{coro}

{\proof}
(i) Clearly, we have that (\ref{argumento}) and (\ref{modulo}) hold. Hence the result follows by Theorem \ref{caracterizacao}.

\noindent (ii)
 If $p$ is odd, using the periodicity of $c_{n}$ and the assumption that $c_{2n}=-c_{2n-1},$ we conclude that $c_{n}=0$, $n \geq 1.$ Hence, $\mu$ is symmetric. Moreover, since $\{m_{n}\}_{n=1}^{\infty}$ is a periodic sequence of period $p$ and $c_{n}=0$, $n \geq 1$, the conditions (\ref{argumento}) and (\ref{modulo}) of Theorem \ref{caracterizacao} can be easily verified. Consequently, the measure $\mu$ has $p\,-$periodic sequence of Verblunsky coefficients.
 \fim\vspace{0.3cm}

The Corollary \ref{periodicaalternante} shows that if we choose the sequence $\{c_{n}\}_{n=1}^{\infty}$ $p\,-$periodic ($p$ even) and  such that $c_{2n}=-c_{2n-1}$, then it is possible, by choosing $\{m_{n}\}_{n=1}^{\infty}$ also $p\,-$periodic, to get a measure $\mu^{(p)}$ whose Verblunsky coefficients are periodic with the same period. Notice that in the case when $c_{2n}=-c_{2n-1}$ and $c_{2n}>0$ (or $c_{2n}<0$) for $n \geq 1$  the sequence $\{c_{n}\}_{n=1}^{\infty}$ has the alternating sign property.

The next theorem provides a geometric characterization for the choice of $\{c_{n}\}_{n=1}^{\infty}$ and $\{m_{n}\}_{n=1}^{\infty}$ considered above.

\begin{theo} \label{teogeometrico}
Let $p$ be an even natural number and $\mu^{(p)}$ be the probability measure associated with the pair $\left\{\{c_{n}\}_{n=1}^{\infty}, \{d_{n}\}_{n=1}^{\infty}\right\}.$ Then, the following statements are equivalent:
\begin{enumerate}
\item[{\rm(i)}] The sequences $\{c_{n}\}_{n=1}^{\infty}$ and $\{m_{n}\}_{n=1}^{\infty}$ are $p\,-$periodic with $c_{2n}=-c_{2n-1},\, n\geq 1.$
\item[{\rm(ii)}] The sequence of Verblunsky coefficients $\{\alpha_{n}\}_{n=0}^{\infty}$ associated with the measure $\mu^{(p)}$ is $p\,-$ periodic. In addition, for $ k \, \in \{0, 1, \ldots, \frac{p-2}{2}\},$ the straight lines connecting $\alpha_{2k}$ to $1$ and $\alpha_{2k+1}$ to $-1$ are parallel.
\end{enumerate}
\end{theo}
{\proof}
(i)$\Rightarrow$(ii) From Corollary \ref{periodicaalternante} it is immediate that $\{\alpha_{n}\}_{n=0}^{\infty}$ is a periodic sequence with period $p.$ On the other hand, by the assumption that $c_{2n}=-c_{2n-1}$ and by (\ref{verblunsky}), for $n \geq 0,$ we have
$$\alpha_{2n}=\frac{b_{2n+1}+ic_{2n+2}}{1+ic_{2n+2}}=1+\lambda_{2n}(-1-ic_{2n+1}),$$
where $\lambda_{2n}=\frac{1-b_{2n+1}}{1+c_{2n+1}^{2}}.$
Similarly, for $n \geq 0$
$$\alpha_{2n+1}=\frac{b_{2n+2}-ic_{2n+2}}{1+ic_{2n+2}}= -1+\lambda_{2n+1}(-1+ic_{2n+2}),$$
where $\lambda_{2n+1}=-\frac{1+b_{2n+2}}{1+c_{2n+2}^{2}}.$

Hence, for each $k \in \{0, 1, \ldots,\frac{p-2}{2}\},$ one can see that $\alpha_{2k} \in r_{2k}$, where $r_{2k}$ is the straight line with parametric equation given by $r_{2k}(t)=1+t(-1-ic_{2k+1}), \, t \in \mathbb{R}.$ 

Similarly, for each $k \in \{0, 1, \ldots,\frac{p-2}{2}\},$ one can see that $\alpha_{2k+1} \in r_{2k+1}$, where $r_{2k+1}$ is the straight line with parametric equation given by $r_{2k+1}(t)=1+t(-1+ic_{2k+2}), \, t \in \mathbb{R}.$

Finally, since $-1-ic_{2k+1}=-1+ic_{2k+2}$ it follows that  $r_{2k}\parallel  r_{2k+1}$, for each $k \in \{0, 1, \ldots,\frac{p-2}{2}\}.$

 \vspace{0.2cm}

 (ii)$\Rightarrow$(i) Let $\alpha_{j}=x_{j}+iy_{j}, \, j=0,1,\ldots, p-1.$ If $j=2k,$ $k \in \{0, 1, \ldots, \frac{p-2}{2}\},$ we can write
 \begin{equation}\label{eq1}
 \alpha_{2k}=1+\lambda_{2k}(-1-i\tilde{c}_{2k+1}),
 \end{equation}
  where
  \begin{equation}\label{eq2}
  \lambda_{2k}=\frac{1-\tilde{b}_{2k+1}}{1+\tilde{c}_{2k+1}^{2}}, \, \  \tilde{c}_{2k+1}=\frac{y_{2k}}{x_{2k}-1}  \quad \mbox{and} \quad \tilde{b}_{2k+1}=1+\frac{(x_{2k}-1)^{2}+y_{2k}^{2}}{x_{2k}-1}.\end{equation}

  Likewise, if $j=2k+1,$ $k \in \{0, 1, \ldots, \frac{p-2}{2}\},$ we can write
 \begin{equation} \label{eq3}
 \alpha_{2k+1}=-1+\lambda_{2k+1}(-1+i\tilde{c}_{2k+2}),
 \end{equation}
  where
  \begin{equation} \label{eq4}
  \lambda_{2k+1}=-\frac{1+\tilde{b}_{2k+2}}{1+\tilde{c}_{2k+2}^{2}}, \, \ \tilde{c}_{2k+2}=\frac{-y_{2k+1}}{1+x_{2k+1}}  \quad \mbox{and} \quad \tilde{b}_{2k+2}=-1+\frac{(1+x_{2k+1})^{2}+y_{2k+1}^{2}}{1+x_{2k+1}}.\end{equation}

  Hence, if we set $\tilde{b}_{n}=1-2\tilde{m}_{n,}$ from $\alpha_{n+p}=\alpha_{n},$ (\ref{eq2}) and (\ref{eq4}) one can see that

  \begin{equation}\label{eq5}
  \tilde{c}_{n+p}=\tilde{c}_{n} \quad \mbox{and} \quad \tilde{m}_{n+p}=\tilde{m}_{n}, \quad n=1, 2, \ldots\, \ .
  \end{equation}

 For each $k  \in \{0, 1, \ldots, \frac{p-2}{2}\},$ let $r_{2k}$ be the straight line connecting $\alpha_{2k}$ to $1$ and $r_{2k+1}$ the straight line connecting $\alpha_{2k+1}$ to $-1.$  Then, from (\ref{eq1}), (\ref{eq3}), (\ref{eq5}) and since $r_{2k}\parallel r_{2k+1},\,\, k  \in \{0, 1, \ldots, \frac{p-2}{2}\},$  it follows that   
  \begin{equation}\label{c2n}\tilde{c}_{2n+2}=-\tilde{c}_{2n+1}, \quad  n=0, 1, \ldots\, \ .
  \end{equation}
 
  Hence, from (\ref{eq1}) to (\ref{c2n}) we have, for $n \geq 1,$
    \begin{equation}\label{eq6}
  \alpha_{2n}=\frac{\tilde{b}_{2n+1}+i\tilde{c}_{2n+2}}{1+i\tilde{c}_{2n+2}} \quad \mbox{and} \quad  \alpha_{2n+1}=\frac{\tilde{b}_{2n+2}-i\tilde{c}_{2n+2}}{1+i\tilde{c}_{2n+2}}.
  \end{equation}
  
  Finally, using the formula (\ref{verblunsky}) for $\alpha_{n}$ and the relation (\ref{eq6}) one can see, by mathematical induction, that for $n \geq 1,$
  $$\tilde{c}_{n}=c_{n} \quad \mbox{and} \quad \tilde{m}_{n}=m_{n}.$$
This completes the proof.
\fim

\vspace{0.3cm}
Observe that Theorem \ref{teogeometrico} shows that to choose a periodic sequence $\{\alpha_{n}\}$ of period $p$ ($p$ even) with $\alpha_{j}$ on certain parallel straight lines is equivalent to choosing the sequences $\{c_{n}\}$ and $\{m_{n}\}$ also $p-$periodic with the additional property $c_{2n+2}=-c_{2n+1}, \, n\geq 0.$ In Fig. 1 and Fig. 2 we show some examples of possible choices for $\{c_{n}\}$ and $\{m_{n}\}.$

\begin{figure}[!ht]
\hspace{0.8cm}\begin{minipage}[t]{0.40\linewidth}
\includegraphics[width=\linewidth]{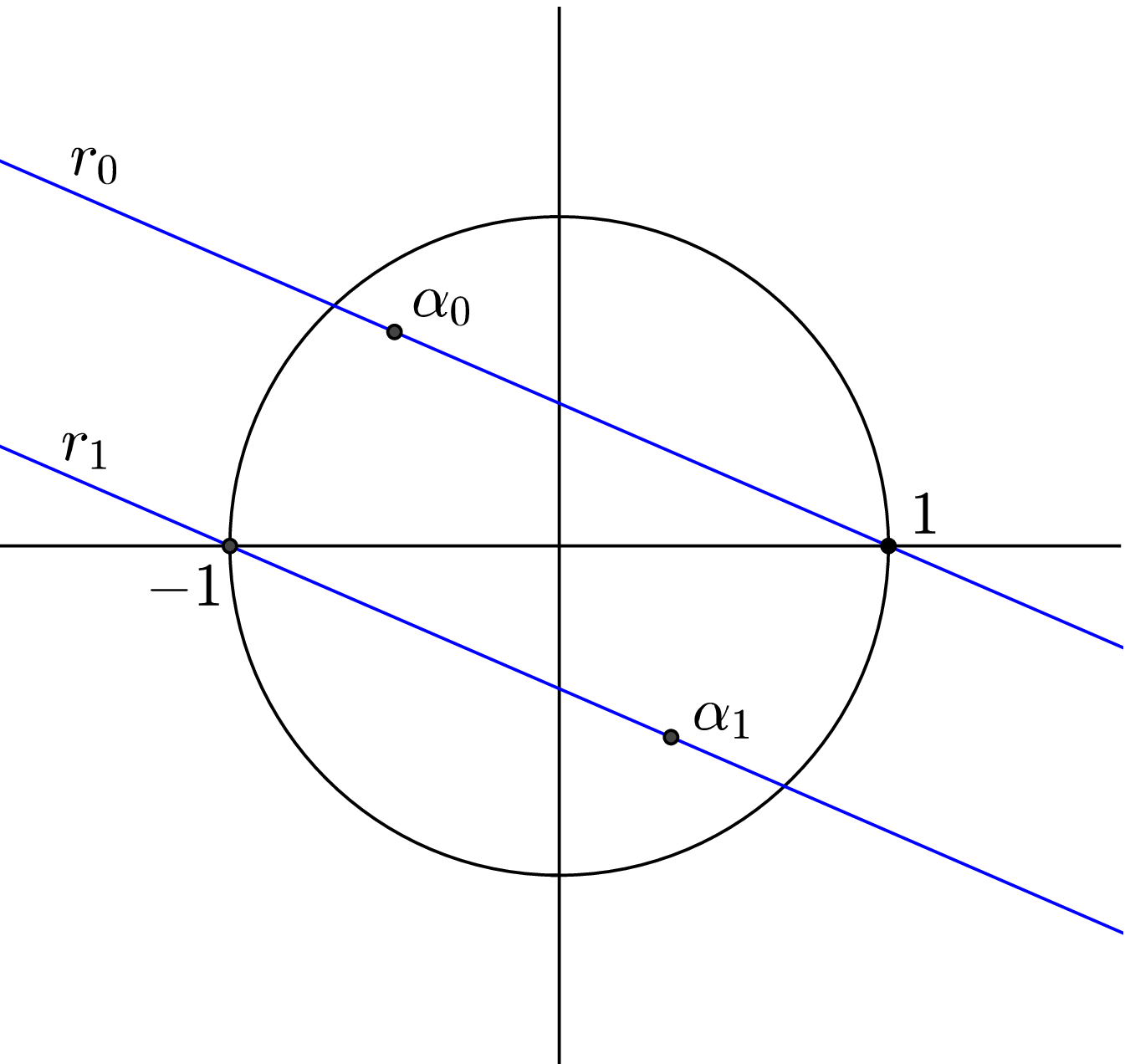}
\caption{{\scriptsize Verblunsky coefficients  associated to the choice $\{c_{n}\}=(-c, c, -c, c, \ldots)$ and $\{b_{n}\}=(b_{1}, b_2, b_1,b_2, \ldots)$, with $c>0.$}}
\label{fig1}
\end{minipage}\hspace{2cm}
\begin{minipage}[t]{0.40\linewidth}
\includegraphics[width=\linewidth]{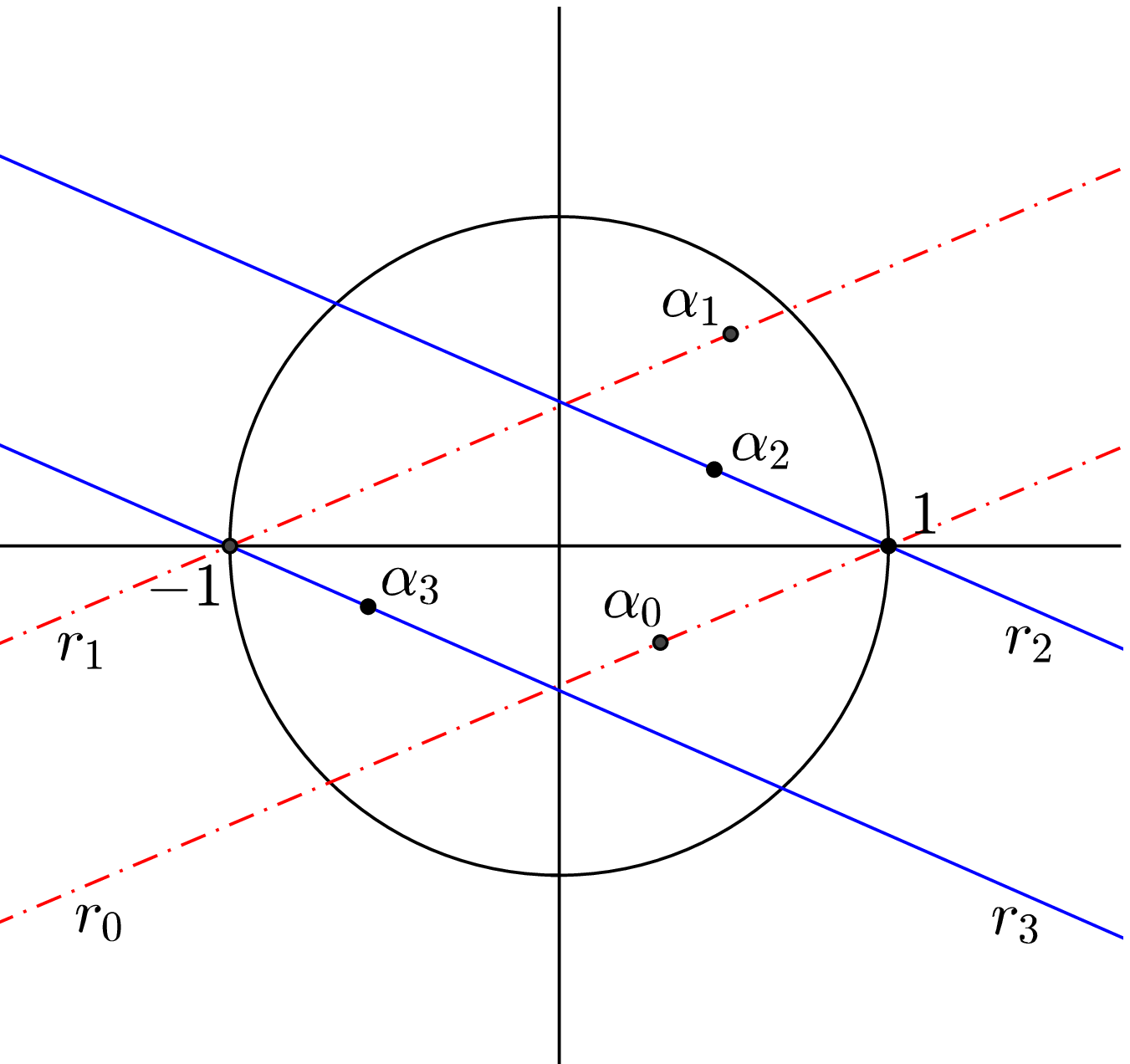}
\caption{{\scriptsize Verblunsky coefficients  associated to the choice $\{c_{n}\}=(-c_2, c_2, -c_4, c_4, -c_2, c_2, \ldots)$ and $\{b_{n}\}=(b_{1}, b_2, b_3, b_{4}, b_{1}, b_2,  \ldots)$, with $c_2 <0$ and $c_4 >0$.}}
\label{fig2}
\end{minipage}
\end{figure}

The following results deal with how to calculate (from the point of view of chain sequences) the pure points and the respective masses of a measure $\mu^{(p)}$ whose associated Verblunsky coefficients are periodic. In \cite{Simon-book-p2} there is another approach to the same problem. 

We begin with a lemma that leads to a characterization of the possible pure points (that we denote by $w$) of the measure $\mu^{(p)}$ in terms of the sequence $\{\tau_{n}(w)\}$ defined in (\ref{taujw}).   
\begin{lema}\label{teotauj}
Let $\mu^{(p)}$ be a probability measure on the unit circle with $p\,-$periodic Verblunsky coefficients. Then, $w$ is a possible pure point of the measure $\mu^{(p)}$ if, and only if, the sequence $\{\tau_{n}(w)\}_{n=0}^{\infty}$ is periodic of period $p.$
\end{lema}

{\proof}
By Theorem \ref{puro1} and Theorem \ref{puro2} we see that $w$ is a possible pure point of $\mu^{(p)}$ if, and only if, $\varphi_{p}(w)-\varphi_{p}^{*}(w)=0.$ Notice that the condition $\varphi_{p}(w)-\varphi_{p}^{*}(w)=0$ is equivalent to $\tau_{p}(w)=1.$

Furthermore, using the periodicity of the sequence $\{\alpha_{n}\}_{n=0}^{\infty}$ and the recurrence relation (\ref{taujw}), we also see that $\tau_{p}(w)=1$ is equivalent to the periodicity of the sequence $\{\tau_{n}(w)\}_{n=0}^{\infty}.$
\fim

\vspace{0.3cm}
The next theorem provides a way to determinate all the pure points of the measure $\mu^{(p)}$ and also, to calculate the mass of each pure point.

\begin{theo}\label{purepoint}
Let $\mu^{(p)}$ be a probability measure on the unit circle with $p\,-$periodic sequence $\{\alpha_{n}\}_{n=0}^{\infty}$ of Verblunsky coefficients. In addition, suppose that $w$ is a point on the unit circle such that $\varphi_{p}(w)-\varphi_{p}^{*}(w)=0.$ Then, $w$ is a pure point of $\mu^{(p)}$ if, and only if,
$$\prod_{j=1}^{p}|1-w\tau_{j-1}(w)\alpha_{j-1}|^{2}<\prod_{j=1}^{p}\left[1-|\alpha_{j-1}|^{2}\right].$$
Moreover, if $w$ is a pure point of  $\mu^{(p)},$ then the mass at this point is given by
$$\mu^{(p)}(\{w\})=\frac{\gamma}{\gamma+\delta}\, ,$$
where $\displaystyle{\delta=\sum_{n=1}^{p}\prod_{j=1}^{n} \frac{|1-w\tau_{j-1}(w)\alpha_{j-1}|^{2}}{1-|\alpha_{j-1}|^{2}}}$
and $\displaystyle{\gamma=1-\prod_{j=1}^{p} \frac{|1-w\tau_{j-1}(w)\alpha_{j-1}|^{2}}{1-|\alpha_{j-1}|^{2}}.}$
\end{theo}

{\proof}
For $j=1, 2, \ldots,$ let $q_{j}=\frac{|1-w\tau_{j-1}(w)\alpha_{j-1}|^{2}}{1-|\alpha_{j-1}|^{2}}.$

By Theorem \ref{teo2} we know that $w$ is a pure point if, and only if, the infinite sum $\lambda(w)=\sum_{n=1}^{\infty}\prod_{j=1}^{n}q_{j}$ is convergent.

\vspace{0.2cm}
By Lemma \ref{teotauj} and by the periodicity of $\{\alpha_{n}\}_{n=0}^{\infty}$ it follows that $q_{j+p}=q_{j}, \, j \geq 1.$

Thus, if $\displaystyle{q=\prod_{j=1}^{p}q_{j},}$ we can write $\lambda(w)$  as
\begin{eqnarray}\label{somatorioq}
\lambda(w)=q_{1}\left(\sum_{n=0}^{\infty}q^{n}\right)+q_{1}q_{2}\left(\sum_{n=0}^{\infty}q^{n}\right)+\cdots+q_{1}q_2\cdots q_{p}\left(\sum_{n=0}^{\infty}q^{n}\right).
\end{eqnarray}

Observe that $\lambda(w)$ is convergent if, and only if, $|q|<1.$ Thus, the first part of the statement follows.

Furthermore, if $|q|<1$ using (\ref{somatorioq}), we have 
\begin{equation}\label{lambdaw}
\lambda(w)=\left(\frac{1}{1-q}\right)\left(\sum_{n=1}^{p}\prod_{j=1}^{n}q_{j}\right)=\frac{\delta}{\gamma}.
\end{equation}

Finally, by Theorem \ref{teo2} and (\ref{lambdaw}), we get
$$\mu^{(p)}(\{w\})=\frac{1}{1+\lambda(w)}=\frac{\gamma}{\gamma+\delta}\,.$$\fim

\vspace{-0.5cm}
\setcounter{equation}{0}
\section{An example}   \label{example}

In this section we  discuss, using the following example, the results obtained in the previous sections.

Let the real sequences $\{c_{n}\}_{n=1}^{\infty}$ and $\{d_{n}\}_{n=1}^{\infty}$ be given by
$$c_{n}=(-1)^nc\quad\mbox{and}\quad d_{n}=(1-m_{n-1})m_n,\quad n\geq 1,$$
where $c\in\mathbb{R}$ and the real sequence $\{m_{n}\}_{n=0}^{\infty}$ is such that $m_{0}=0$,
$$m_{2n-1}=\frac{1-b_1}{2}\quad\mbox{and}\quad m_{2n}=\frac{1-b_2}{2},\quad n\geq1,$$
with $b_1, b_2\in\mathbb{R}$ and $|b_1|, |b_2|<1$.

Notice that, if $c \neq 0,$  $\{c_{n}\}_{n=1}^{\infty}$ has the alternating sign property and that $\{d_{n}\}_{n=1}^{\infty}$ is a positive chain sequence, with $\{m_{n}\}_{n=0}^{\infty}$ being its minimal parameter sequence. Moreover, $\{c_{n}\}_{n=1}^{\infty}$ and $\{m_{n}\}_{n=1}^{\infty}$ are periodic sequences of period 2.

By Theorem \ref{teo1}, associated with the pair $\left\{\{c_{n}\}_{n=1}^{\infty}, \{d_{n}\}_{n=1}^{\infty}\right\},$ there exists an unique probability measure, say $\mu^{(2)},$ on the unit circle. Furthermore, from Corollary \ref{periodicaalternante} follows that the sequence of Verblunsky coefficients of $\mu^{(2)}$ is periodic with period $2$ (in Fig. 1, it is illustrated the position of these coefficients for the case $c>0$).

From the definition of $\{c_{n}\}_{n=1}^{\infty}$ one can also see that
\begin{equation}\label{tau1}
\tau_{2n}=1 \quad \mbox{and} \quad \tau_{2n+1}=\frac{1+ic}{1-ic}\,\,, \ \ n\geq 0.
\end{equation}
Thus, from (\ref{verblunsky}) we have, for $n\geq 0$,
$$\alpha_{2n}=\frac{b_{1}+ic}{1+ic}=\frac{(b_{1}+c^{2})+ic(1-b_{1})}{1+c^{2}} \quad \mbox{and} \quad
\alpha_{2n+1}=\frac{b_{2}-ic}{1+ic}=\frac{(b_{2}-c^{2})-ic(1+b_{2})}{1+c^{2}} \, \,.$$

In this case, since $p=2,$ we have $\Delta(z)=z^{-1}{\rm Tr}(T_{2}(z)).$ By (\ref{transferencia}) and (\ref{matrizA})
\begin{eqnarray*}T_{2}(z)=(1-|\alpha_{0}|^{2})^{-1/2}(1-|\alpha_{1}|^{2})^{-1/2}\left(\begin{array}{cc}
                        z &  -\overline{\alpha}_{1}\\
                         -\alpha_{1} z & 1
           \end{array}\right)\left(\begin{array}{cc}
                        z &  -\overline{\alpha}_{0}\\
                         -\alpha_{0} z & 1
           \end{array}\right).\end{eqnarray*}

Hence, computing $\Delta(e^{i\theta})$ one can see that, for $\theta \in [0,2\pi),$
$$\Delta(e^{i\theta})=2\left\{\frac{1+c^{2}}{[(1-b_{1}^{2})(1-b_{2}^{2})]^{1/2}}\cos \theta+\frac{b_{1}b_{2}-c^{2}}{[(1-b_{1}^{2})(1-b_{2}^{2})]^{1/2}}\right\} $$
and, consequently
$$\sqrt{4-\Delta^{2}(e^{i\theta})}=2\sqrt{1-\left[\frac{(1+c^{2})\cos \theta+b_{1}b_{2}-c^{2}}{(1-b_{1}^{2})^{1/2}(1-b_{2}^{2})^{1/2}}\right]^{2}}\,\,.$$

Furthermore, considering the normalized orthogonal polynomials $\varphi_{2}(z)=\kappa_{2}\phi_{2}(z)$ one can also verify that
{\small $$\varphi_{2}(z)=\frac{1}{(1-b_{1}^{2})^{1/2}(1-b_{2}^{2})^{1/2}}\left\{(1+c^{2})z^{2}+[(b_{1}b_{2}-b_{1}-2c^{2})+ic(b_{2}+1)]z+[(c^{2}-b_{2})-ic(b_{2}+1)]\right\}$$}
and, consequently for $\theta \in [0,2\pi),$
$${\rm Im}(e^{-i\theta}\varphi_{2}(e^{i\theta}))=\frac{(1+b_{2})[\sin\theta+c(1-\cos\theta)]}{(1-b_{1}^{2})^{1/2}(1-b_{2}^{2})^{1/2}}\,.$$

Hence, from Theorem \ref{abscont}, the weight function $w(\theta)$ associated to $\mu^{(2)}$ is such that

$$w(\theta)=\frac{\sqrt{(1-b_{1}^{2})(1-b_{2}^{2})-[(1+c^{2})\cos\theta+b_{1}b_{2}-c^{2}]^{2}}}{|(1+b_{2})[\sin\theta+c(1-\cos\theta)]|}\,.$$

Now we need to compute the bands $B_{1}$ and $B_{2}$ for the measure $\mu^{(2)}.$ By solving the equation $\Delta(e^{i\theta})=2$ we find the solutions
$$\theta_{1}^{+}=\arccos\left(\frac{(1-b_{1}^{2})^{1/2}(1-b_{2}^{2})^{1/2}+c^{2}-b_{1}b_{2}}{1+c^{2}}\right) \quad \mbox{and} \quad  \theta_{2}^{+}=2\pi-\theta_{1}^{+}.$$
Likewise, by solving $\Delta(e^{i\theta})=-2$ we find
$$\theta_{1}^{-}=\arccos\left(\frac{c^{2}-(1-b_{1}^{2})^{1/2}(1-b_{2}^{2})^{1/2}-b_{1}b_{2}}{1+c^{2}}\right) \quad \mbox{and} \quad  \theta_{2}^{-}=2\pi-\theta_{1}^{-}.$$
Thus, each band $B_{j}$ is determined by the points $z_{j}^{+}=e^{i\theta_{j}^{+}}$ and $z_{j}^{-}=e^{i\theta_{j}^{-}}, \, j \in \{1,2\}.$

To determine the possible pure points of $\mu^{(2)},$ by Theorem \ref{puro1} and Theorem \ref{puro2} we need to solve the equation $\varphi_{2}(z)-\varphi_{2}^{*}(z)=0,$ whose solutions are $w_{1}=1$ and $w_{2}=\frac{c^{2}-1}{1+c^{2}}-i\frac{2c}{1+c^{2}}\,.$

Now looking at the bands $B_{j}$ and for the possible pure points $w_{j}$, it is not hard to see that the measure $\mu^{(2)}$ is always supported on  $\mathcal{C}_{1} \cup \mathcal{C}_{2},$ in accordance with Theorem \ref{teo-geral-alternantes}. 

Finally, we give a complete characterization about the singular part of the measure $\mu^{(2)}$ in terms of the parameters $b_{1},$ $b_{2}$ and $c.$

Firstly, we analyze the point $w_{1}=1.$ Notice that $\tau_{n}(w_{1})=\tau_{n}$ given by (\ref{tau1}) is periodic of period $2,$ according to Lemma \ref{teotauj}. From Theorem \ref{purepoint} one can see that $w_{1}$ is a pure point of $\mu^{(2)}$ if, and only if,  $\,b_{1}+b_{2}>0.$ Moreover, if $\displaystyle{\delta_{1}=\sum_{n=1}^{2}\prod_{j=1}^{n} \frac{|1-\tau_{j-1}\alpha_{j-1}|^{2}}{1-|\alpha_{j-1}|^{2}}}$ and $\displaystyle{\gamma_{1}=1-\prod_{j=1}^{2} \frac{|1-\tau_{j-1}\alpha_{j-1}|^{2}}{1-|\alpha_{j-1}|^{2}},}$ again by Theorem \ref{purepoint} we obtain $$\mu^{(2)}(\{w_{1}\})=\frac{\gamma_{1}}{\gamma_{1}+\delta_{1}}=\frac{b_{1}+b_{2}}{1+b_{2}}\, .$$

Now we consider the point $w_{2}=\frac{c^{2}-1}{1+c^{2}}-i\frac{2c}{1+c^{2}}=-\frac{1+ic}{1-ic}.$  From Corollary \ref{rotacao}, if $\tilde{\mu}(z)=\mu(w_{2}z),$ we have $\tilde{c}_{n}=c, \, n\geq 1.$ Moreover,
$$\tilde{\tau}_{n}=\prod_{k=1}^{n}\frac{1-i\tilde{c}_{k}}{1+i\tilde{c}_{k}}=\left(\frac{1-ic}{1+ic}\right)^{n}, \quad n\geq 1.$$

On the other hand, it is known (see, for example, \cite{Costa-Felix-Ranga-2013}) that $\tilde{\tau}_{n}=w_{2}^{-n}\tau_{n}(w_{2}),\, n \geq 0.$ Hence, one can see that $\tau_{n}(w_{2})=(-1)^{n}, \, n \geq 0.$ Thus, it follows that $\tau_{n}(w_{2})$ is periodic of period $2,$ according to Lemma \ref{teotauj}.

From Theorem \ref{purepoint},  $w_{2}$ is a pure point of $\mu^{(2)}$ if, and only if, $b_{2}-b_{1}>0.$ Moreover, if $\displaystyle{\delta_{2}=\sum_{n=1}^{2}\prod_{j=1}^{n} \frac{|1-w_{2}\tau_{j-1}(w_{2})\alpha_{j-1}|^{2}}{1-|\alpha_{j-1}|^{2}}}$ and $\displaystyle{\gamma_{2}=1-\prod_{j=1}^{2} \frac{|1-w_{2}\tau_{j-1}(w_{2})\alpha_{j-1}|^{2}}{1-|\alpha_{j-1}|^{2}},}$ we obtain
 $$\mu^{(2)}(\{w_{2}\})=\frac{\gamma_{2}}{\gamma_{2}+\delta_{2}}=\frac{b_{2}-b_{1}}{1+b_{2}}\, .$$
 
Finally, some possible situations are illustrated in the figures $3, 4, 5$ and $6$.  
\vspace{0.3cm}
\begin{figure}[h]
\hspace{0.8cm}\begin{minipage}[t]{0.40\linewidth}
\includegraphics[width=\linewidth]{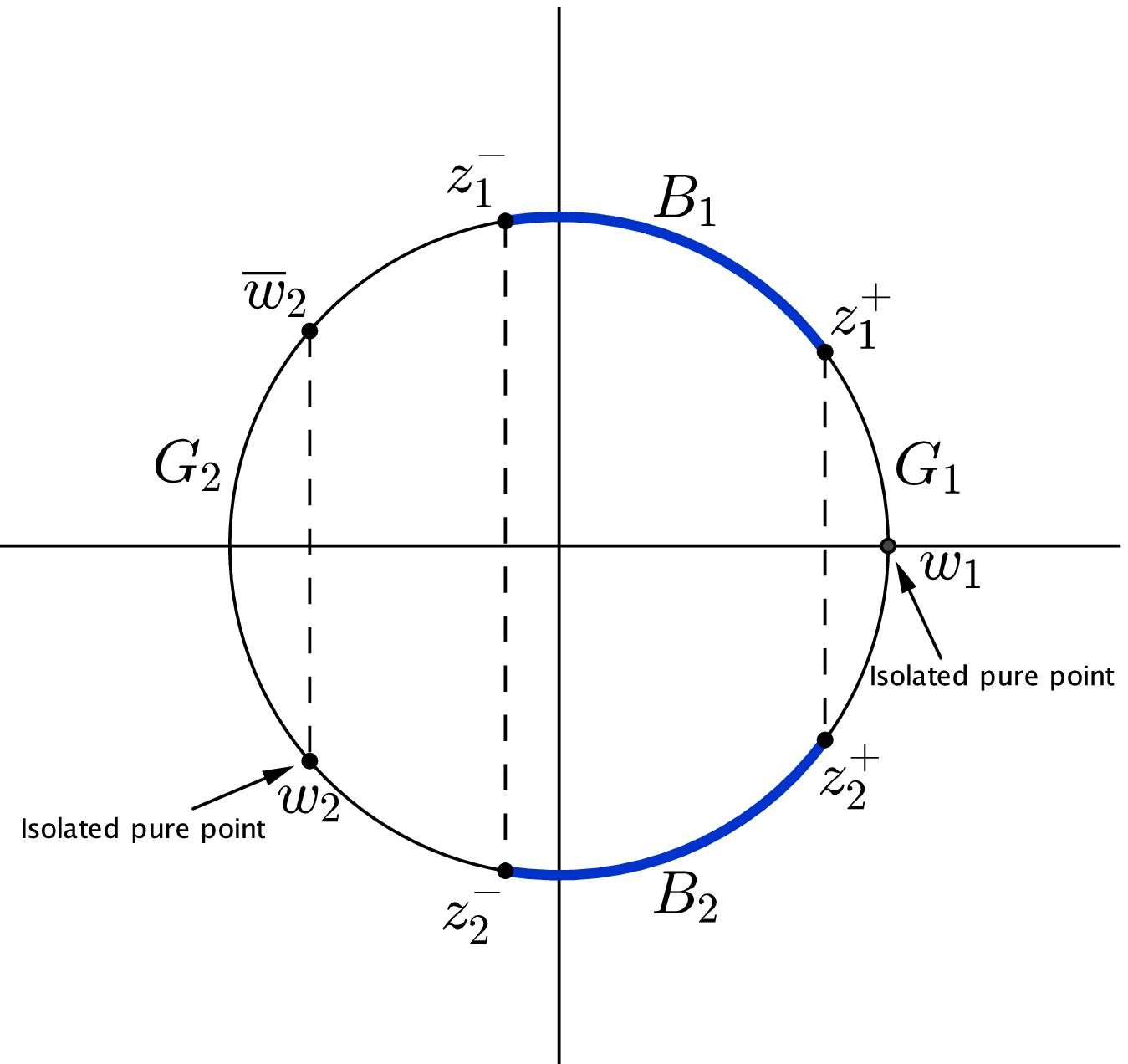}
\caption{{\scriptsize Support of $\mu^{(2)}$ in the case $0<c<1$ and $b_2>b_1>0$.}}
\label{fig3}
\end{minipage}\hspace{2cm}
\begin{minipage}[t]{0.40\linewidth}
\includegraphics[width=\linewidth]{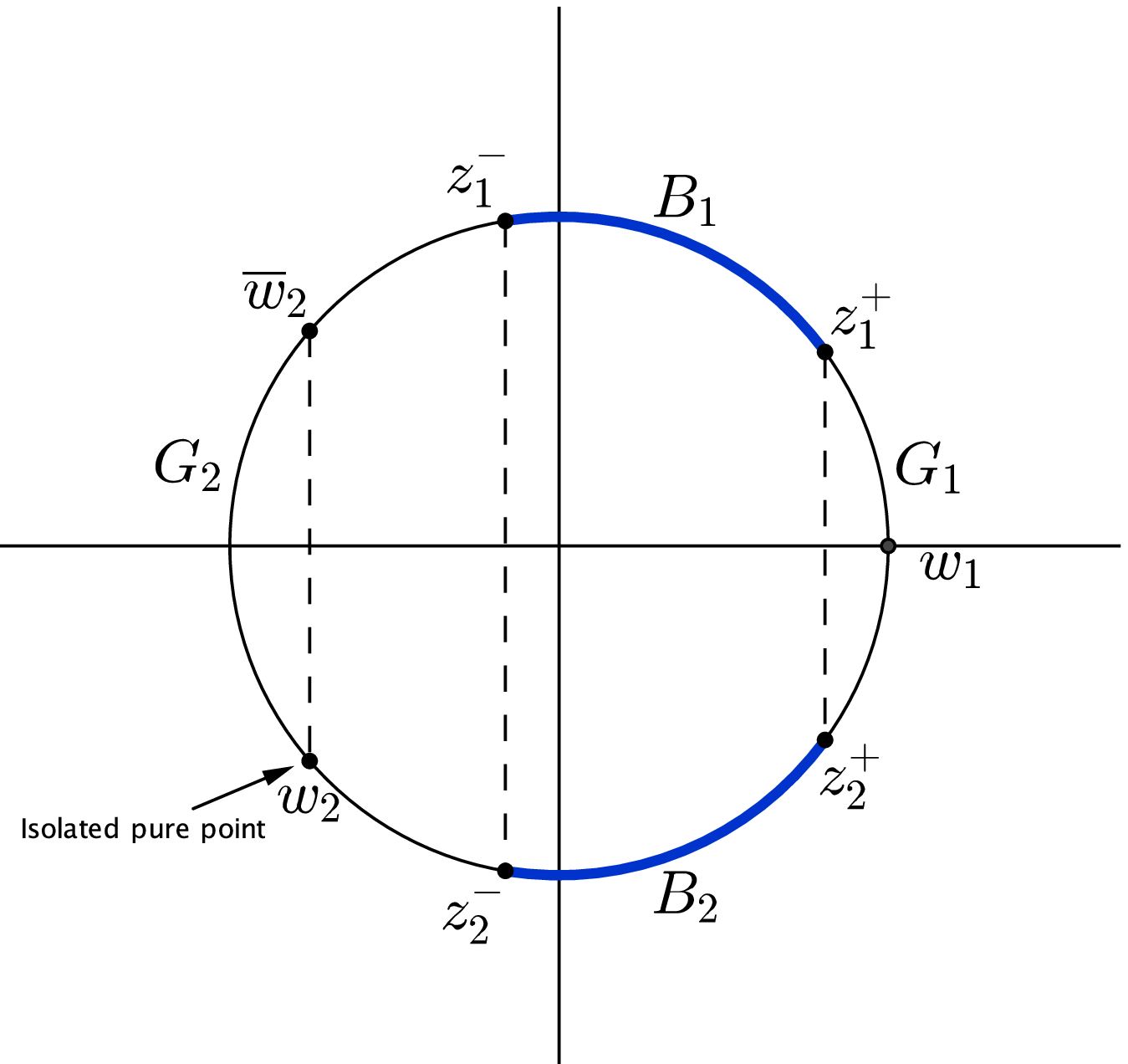}
\caption{{\scriptsize Support of $\mu^{(2)}$ in the case $0<c<1$ and $0<b_2 \leq -b_1$ .}}
\label{fig4}
\end{minipage}
\end{figure}

\begin{figure}[h]
\hspace{0.8cm}\begin{minipage}[t]{0.40\linewidth}
\includegraphics[width=\linewidth]{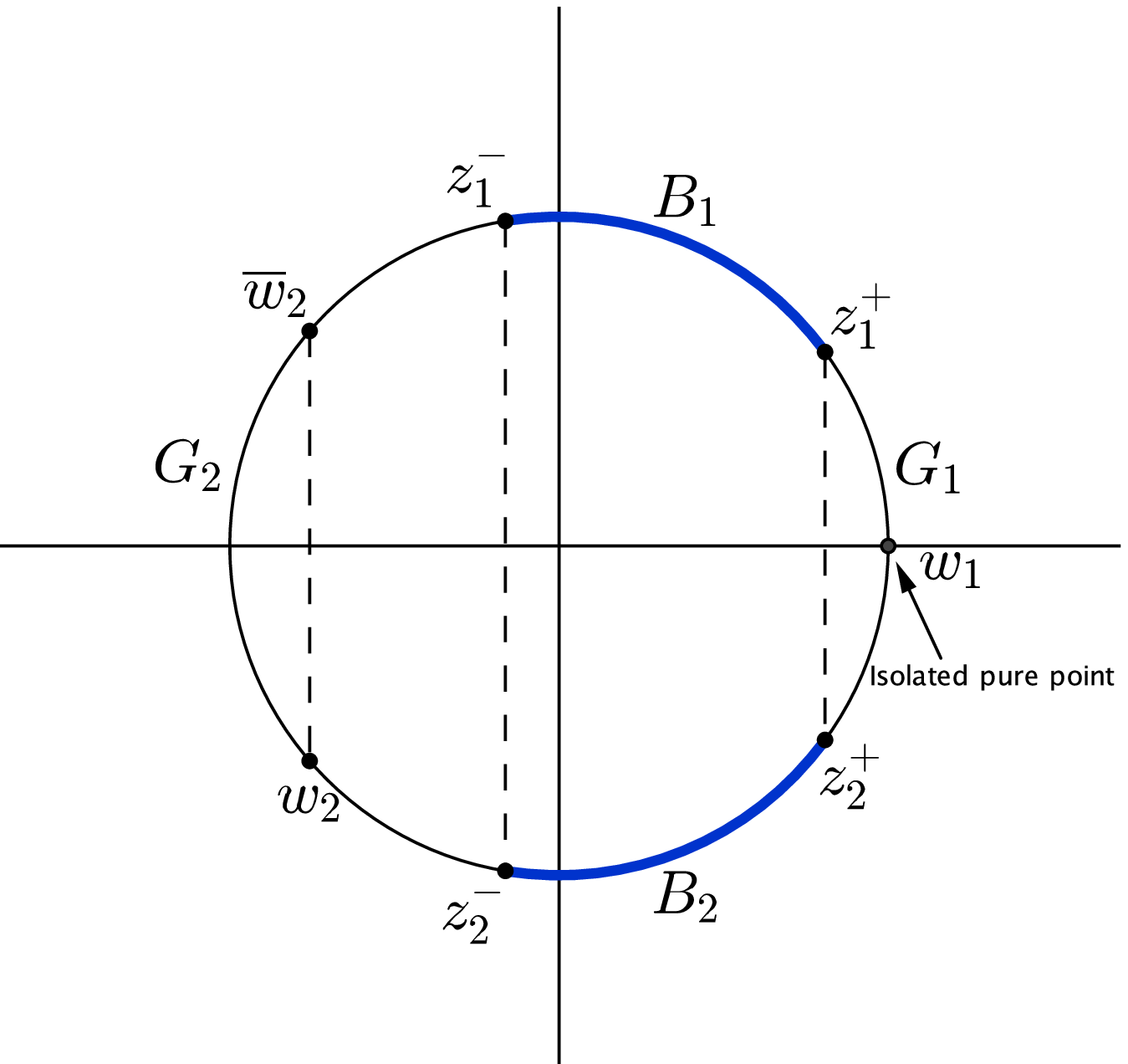}
\caption{{\scriptsize Support of $\mu^{(2)}$ in the case $0<c<1$ and $0<b_2\leq b_1$.}}
\label{fig5}
\end{minipage}\hspace{2cm}
\begin{minipage}[t]{0.40\linewidth}
\includegraphics[width=\linewidth]{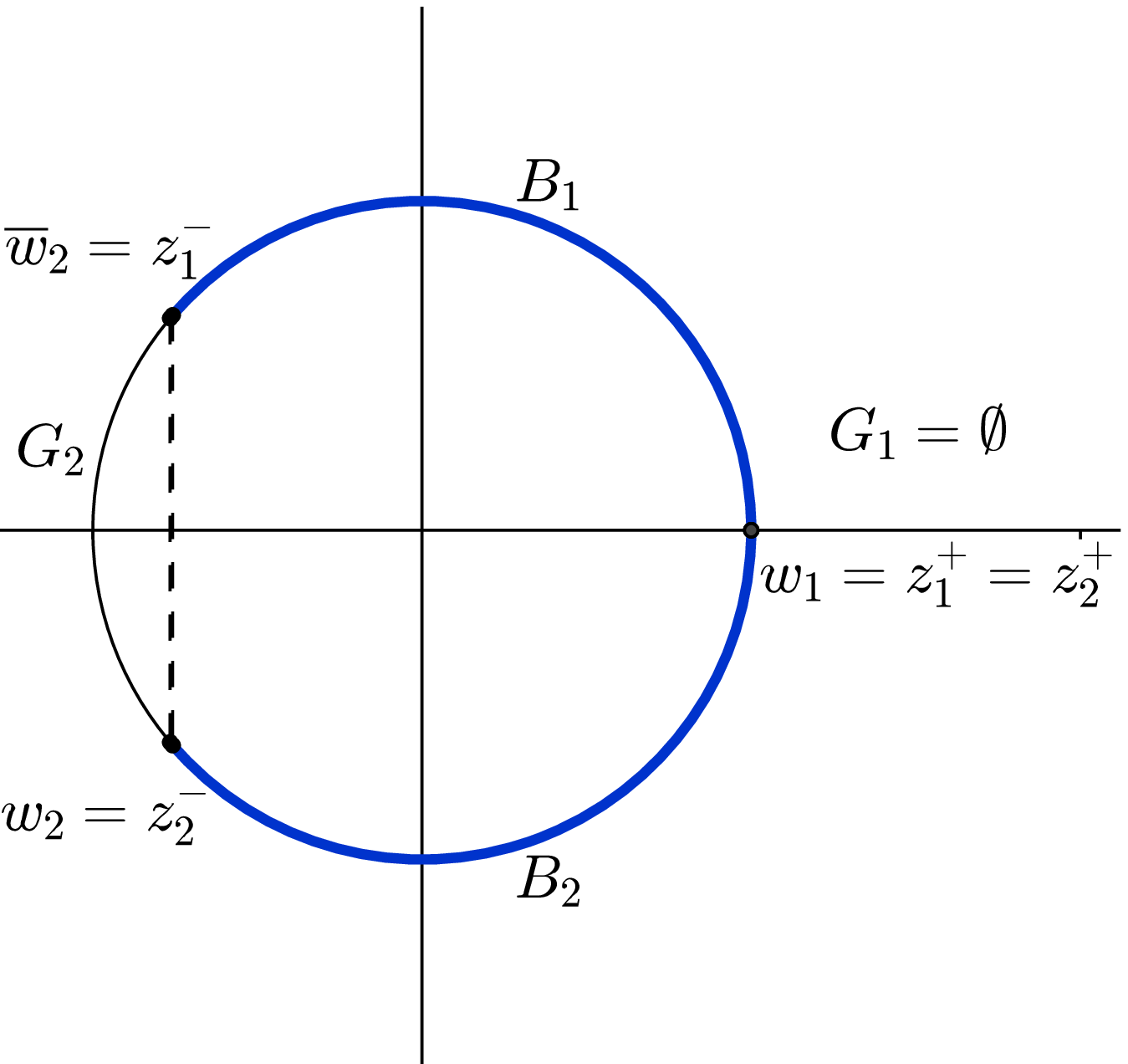}
\caption{{\scriptsize Support of $\mu^{(2)}$ in the case $0<c<1$ and $b_1=b_2=0$ .}}
\label{fig6}
\end{minipage}
\end{figure}
\newpage



\begin{thebibliography}{99}

\bibitem{BracMcCabPerezRanga-MCOM2015} C. F. Bracciali, J.H. McCabe, T.E. Per\'{e}z and A. Sri Ranga, A class of orthogonal functions given by a three term recurrence formula, {\em Math. Comp.}, (to appear).

\bibitem{BreuerSimon-2010}  J. Breuer, E. Ryckman and B. Simon,  Equality of the spectral and dynamical definitions of reflection, {\em Comm. Math. Phys.}, 295 (2010), 531-550.

\bibitem{Castillo-Costa-Ranga-Veronese-2014} K. Castillo, M. S. Costa, A. Sri Ranga and D. O. Veronese, A Favard type theorem for orthogonal polynomials on the unit circle from a three term recurrence formula, {\em J. Approx. Theory}, 184 (2014), 146-162.

\bibitem{Castil-Garza-Marcell-2011} K. Castillo, L. Garza and F. Marcell\'{a}n, Perturbations on the subdiagonals of Toeplitz matrices, {\em  Linear Algebra Appl.}, 434 (2011), 1563-1579.

\bibitem{Chihara-book} T. S. Chihara,  {``An Introduction to Orthogonal Polynomials''}, Mathematics and its Applications Series,  Gordon and Breach, New York, 1978.

\bibitem{Costa-Felix-Ranga-2013} M. S. Costa, H. M. Felix and A. Sri Ranga, Orthogonal polynomials on the unit circle and chain sequences, {\em J. Approx. Theory}, 173 (2013), 14-32.

\bibitem{Costa-Ranga-Godoy-2011} M. S. Costa, E. Godoy, R. L. Lamblém and A. Sri Ranga, Basic hypergeometric functions and orthogonal Laurent polynomials, {\em Proc. Amer. Math. Soc.}, 140 (2011), 2075-2089.

\bibitem{DimRan-2013} D. K. Dimitrov and A. Sri Ranga, Zeros of a family of hypergeometric  para-orthogonal polynomials on the unit circle, {\em Math. Nachr.}, 286 (2013), 1778-1791.

\bibitem{Nevai} T. Erdélyi, P. Nevai, J. Zhang and J. Geronimo,  A simple proof of {``Favard's theorem''} on the unit circle, {\em Atti Semin. Mat. Fis. Univ. Modena Reggio Emilia},  39 (1991), 551-556. Also in {``Trends in functional analysis and approximation theory''} (Acquafredda di Maratea, 1989), 41-46, Univ. Moderna Reggio Emilia, Moderna, 1991.

\bibitem{Geronimus-1941} Ya. L. Geronimus, On polynomials orthogonal on the unit circle, on trigonometric moment problem, and on allied Carathéodory and Schur functions, {\em Mat. Sb.}, 15 (1944), 99-130.

\bibitem{Ismail-book} M. E. H. Ismail, {``Classical and Quantum Orthogonal Polynomials in One Variable''}, Encyclopedia of Mathematics and its Applications Vol. 98, Cambridge Univ. Press, Cambridge,  2005.

\bibitem{Kheif-Golins-Pehers-Yudit-2011} A. Kheifets, L. Golinskii, F. Peherstorfer and P. Yuditskii, Scattering theory for CMV Matrices: uniqueness, Helson-Szeg\H{o} and strong Szeg\H{o} theorems, {\em Integral Equations Operator Theory}, 69 (2011), 479-508.

\bibitem{Pehers-2011} F. Peherstorfer, Positive trigonometric quadrature formulas and quadrature on the unit circle,  {\em  Math. Comp.}, 80 (2011), 1685-1701.

\bibitem{Pehers-Stein-1996} F. Peherstorfer and R. Steinbauer, Orthogonal polynomials on arcs of the unit circle, I, {\em J. Approx. Theory}, 85 (1996), 140-184.

\bibitem{Pehers-Stein-1997} F. Peherstorfer and R. Steinbauer, Orthogonal polynomials on arcs of the unit circle, II. Orthogonal polynomials with periodic reflection coefficients, {\em J. Approx. Theory}, 88 (1997), 316-353.

\bibitem{Pehers-Volb-Yudit-2011} F. Peherstorfer, A. Volberg and P. Yuditskii, CMV matrices with asymptotically constant coefficients. Szeg\H{o}-Blaschke class, scattering theory, {\em J. Funct. Anal.}, 256 (2009), 2157-2210.

\bibitem{Simanek-2012} B. Simanek, Week convergence of CD kernels: a new approach on the circle and real line, {\em J. Approx. Theory}, 164 (2012), 204-209.

\bibitem{Simon-book-p1} B. Simon, {``Orthogonal Polynomials on the Unit Circle. Part 1. Classical Theory''},  Amer. Math. Soc. Colloq. Publ., vol. 54, part 1, Amer. Math. Soc., Providence, RI,  2005.

\bibitem{Simon-book-p2} B. Simon, {``Orthogonal Polynomials on the Unit Circle. Part 2. Spectral Theory''},  Amer. Math. Soc. Colloq. Publ., vol. 54, part 2, Amer. Math. Soc., Providence, RI,  2005.

\bibitem{Szego-book-1939} G. Szeg\H{o}, {``Orthogonal Polynomials''} Amer. Math. Soc. Colloq. Publ., vol. 23,  Amer. Math. Soc., Providence, RI,  1975. Fourth Edition.

\bibitem{Tsujimoto-Zhedanov-2009} S. Tsujimoto and A. Zhedanov, Elliptic hypergeometric Laurent biorthogonal polynomials with a dense point spectrum on the unit circle, {\em SIGMA Symmetry Integrability Geom. Methods Appl.}, 5 (2009), 30p.


\bibitem{Wall} H. S. Wall, {``Analytic Theory of Continued Fractions''}, D. van Nostrand, 1948.

\end{thebibliography}
\end{document}